\documentclass[amsmath,english,a4paper,graphicx,12pt]{article}
\usepackage{doc}
\usepackage{graphicx}
\usepackage{amsmath,amsthm}
\usepackage{amsfonts}
\usepackage{amssymb}
\usepackage{babel}
\usepackage{latexsym}
\usepackage{mathrsfs}
\usepackage{color}

\newtheorem{prop}{\mbox{~}\hskip.3cm Proposition}
\newtheorem{lemma}{\mbox{~}\hskip.3cm Lemma}
\newtheorem{Theorem}{\mbox{~}\hskip.3cm Theorem}
\newtheorem{Definition}{\mbox{~}\hskip.3cm Definition}
\newtheorem{cor}{\mbox{~}\hskip.3cm Corollary}
\newtheorem{Remark}{\mbox{~}\hskip.3cm Remark}
\newtheorem{Example}{\mbox{~}\hskip.3cm Example}
\newcounter{tmp}
\DeclareMathOperator{\dist}{dist}

\begin{document}

\title{Quadrature formulae for the positive real axis in the setting of Mellin analysis:
Sharp error estimates in terms of the Mellin distance 
}


\author{Carlo Bardaro, \thanks{
Department of Mathematics and Computer Sciences, University of Perugia,
via Vanvitelli 1, I-06123 Perugia, Italy, e-mail: 
carlo.bardaro@unipg.it} \and
Paul L. Butzer, \thanks{Lehrstuhl A fuer Mathematik, RWTH Aachen, Templergraben 55, Aachen, D-52056, Germany, e-mail: butzer@rwth-aachen.de}\and
Ilaria Mantellini, \thanks{Department of Mathematics and Computer Sciences, University of Perugia,
via Vanvitelli 1, I-06123 Perugia, Italy, e-mail: 
mantell@dmi.unipg.it}\and Gerhard Schmeisser\thanks{Department Mathematik, FAU Erlangen-N\"{u}rnberg, Cauerstr. 11, 91058 Erlangen, Germany, 
email: schmeisser@mi.uni-erlangen.de}
}

\date{}

\maketitle
\noindent
\vskip0,3cm

{\bf Abstract}.
The general Poisson summation formula of Mellin analysis can be considered as a quadrature formula for  the positive real axis with remainder. 
For Mellin bandlimited  functions it becomes an exact quadrature formula. Our main aim  
is to study the  speed of convergence to zero of the remainder for a function $f$ in terms of its distance from a space of Mellin bandlimited functions. 
The resulting estimates turn out to be of best possible order. Moreover, we characterize certain rates of convergence in terms of Mellin--Sobolev
and Mellin--Hardy type spaces that contain $f$. Some numerical experiments illustrate and confirm these results.
\noindent
\vskip0,3cm
{\bf Keywords}: Mellin transforms, Mellin--Poisson summation formulae,  polar-analytic functions, quadrature formulae  
\vskip0,3cm
\noindent
{\bf AMS Subject Classification.}
 41A80, 65D30, 65D32,  26A33.

\section{Introduction}
\label{intro}
The general Poisson summation formula (see, e.g., \cite{BN}, \cite{BDFHSS}) is a fundamental mathematical tool, being interconnected with pivotal theorems of 
Fourier analysis, signal theory, numerical analysis and number theory.
In particular, it is an important formula in the study of shift-invariant spaces (see \cite{BDR} and \cite[Chap.~13, Sect.~7]{DL}).
Furthermore it can be interpreted as a quadrature formula with remainder over the whole 
real line. In fact, it is the compound trapezoidal rule on $\mathbb{R}$. 
When appropriately scaled, it is exact for bandlimited functions for which it has features of a Gaussian quadrature formula  
(see \cite{RS}, \cite{SCH72}, \cite{SCH77}, \cite{SCH0}, \cite[\S\,2.11.2]{BSS-Marvasti}).  Recently precise estimates of the remainder in terms of a new metric
have been established (see \cite{BSS} which continues  earlier results in \cite{DRS}). 

For functions $f$ defined on the positive real axis the following elegant formula 
$$\int_0^\infty f(x)dx = \frac{1}{\sigma}\sum_{k=-\infty}^\infty f(e^{k/\sigma}) e^{k/\sigma} + R_\sigma[f] \qquad (\sigma>0)$$
was proposed in \cite{S}. It was deduced from the compound trapezoidal rule on $\mathbb{R}$ by a conformal transformation that maps a strip 
symmetrical to the real line of the complex plane
onto a sector symmetrical to the positive real axis. Furthermore, estimates for the remainder
 $R_\sigma[f]$ were obtained in the case when $f$ extends to an analytic function on that sector. Employing Mellin transform arguments, the above formula 
was  studied in depth and generalized in \cite{SCH}. This formula has features of a Gaussian quadrature formula on $\mathbb{R}^+$, in which the role of
polynomials is taken by Mellin bandlimited functions (see \cite[Sec.\ II--III]{SCH}). For some interesting contributions to quadrature on 
$\mathbb{R}^+$ involving the classical 
Gaussian quadrature
formulas, namely the ones with respect to polynomials, see \cite{Gautschi}, \cite{MM}, \cite{MNM}.  

In the present paper we extend the Fourier results of \cite{BSS} to the case of functions belonging to the space 
$X_c:=\{f:\mathbb{R}^+ \rightarrow \mathbb{C}: (\cdot)^{c-1}f(\cdot) \in L^1(\mathbb{R}^+)\},$ $c$ being a fixed real constant. 
Our approximate quadrature formula is then of the form
\begin{eqnarray}\label{1}
\int_0^\infty f(x) x^{c-1}dx = \frac{1}{\sigma}\sum_{k=-\infty}^\infty f(e^{k/\sigma}) e^{kc/\sigma} + R_{c,\sigma}[f]
\end{eqnarray}
with $\sigma >0.$ Now the starting point is the Mellin--Poisson summation formula, the Mellin counterpart of the Poisson summation formula, 
first stated in \cite{BJ2}.

Just as in Fourier analysis, the Mellin--Poisson summation formula is an exact quadrature formula for Mellin bandlimited functions. 
We study estimates of $R_{c,\sigma}[f]$ in Mellin--Sobolev spaces, also of fractional type, and in certain Hardy type spaces, both in terms of the 
distance functional $\dist_\infty(f, B^1_{c,\sigma})$ introduced in \cite{BBMS2}. The final results give some interesting equivalence theorems 
which characterize function spaces in terms of the speed of convergence to zero of the remainders in (\ref{1}). This generalizes and completes the 
results given in \cite{SCH} by employing a modified notion of analyticity introduced in \cite{BBMS4}, called ``polar-analyticity'', and a generalization 
of the Paley--Wiener theorem in the Mellin frame, also established in \cite{BBMS4} (also see \cite{BBMS} and \cite{BBMS3}).\footnote{%
The definition of polar analyticity of a complex-valued function $f(r,\theta)$ on a domain in $\mathbb{R}^+\times\mathbb{R}$,   
presented in Definition~\ref{def1}, arises
by treating the polar coordinates as cartesian coordinates, and it leads to the classical Cauchy-Riemann 
equations when written in their polar form.  Although others may have 
come across this concept too,  it does not seem to have been used in a 
systematic way as yet.  In Mellin analysis it turns out to be very 
helpful as an efficient approach which is independent of Fourier 
analysis. In particular, it leads to a precise and simple analysis for 
functions defined over the complex logarithm, via the helicoidal surface.
The counterpart of the anchor-theorem of complex analysis, Cauchy's integral formula,
 has already been established for polar-analytic functions, see \cite[Proposition~3.3]{BBMS3}. It
 would indeed be a worthwhile project to develop a complete, independent complex
 analysis in which polar analyticity plays the central role of classical
 analyticity.}

The characterizations are obtained by two different approaches (see Section 4): The first one applies to the so-called Mellin-even part of $f$ 
and makes use of the M\"{o}bius  inversion formula---a tool of number theory first used in numerical analysis by Lyness \cite{LY}, Brass \cite{BRA} 
and Loxton--Sanders \cite{LS}, \cite{LS1} for 
deducing properties of a function $f$ from its remainders in a quadrature formula. The second approach is based on a Mellin--Parseval formula
(see \cite[Theorem~3.2, formula (45)]{BJ4}) and leads to equivalence 
theorems which involve the translated functions $\tau^c_hf$ in place of $f$. Here $\tau^c_h$ is the Mellin translation operator. 

Section~2 has been designed for the reader's convenience. We recall basic notions of Mellin analysis and some earlier results that will be employed
later on. In Section~3, we establish estimates of the remainder $R_{c, \sigma}[f]$ in Mellin--Sobolev spaces of both, integer and fractional order.
Section~4 is devoted to the announced equivalence theorems.
In Section 5, we illustrate the previous results by applying them to several novel examples. The final section is devoted to a short biography of the late 
Professor Helmut Brass who significantly promoted the theory of quadrature by inspiring research articles, books and Oberwolfach conferences conducted
by him.

\section{ Basic notions and preliminary results}
\label{sec:1}
In this section we present some basic definitions and preliminary results concerning the Mellin transform and the Poisson summation formula.

\subsection{The Mellin transform and related concepts}
\label{sec:2.1}
Let $C(\mathbb{R}^+)$ be the space of all continuous complex-valued functions defined on $\mathbb{R}^+,$ and $C^{(r)}(\mathbb{R}^+)$ be the space of 
all functions 
in $C(\mathbb{R}^+)$ with a derivative of order $r$ in $C(\mathbb{R}^+).$ Analogously, we write $C^\infty(\mathbb{R}^+)$ for the space of all 
infinitely differentiable functions.
By $L^1_{\rm{loc}}(\mathbb{R}^+)$, we denote the space of all measurable functions which are integrable on every bounded subinterval of $\mathbb{R}^+.$

For $1\leq p < +\infty,$ let $L^p(\mathbb{R}^+)$ be the space of all Lebesgue measurable and $p$-integrable complex-valued functions defined on 
$\mathbb{R}^+$. This space is endowed with the usual norm $\|\cdot\|_p.$ Analogous notations are used for functions 
defined on $\mathbb{R}.$

For $c \in \mathbb{R},$ we introduce the space
$$X_c = \{ f: \mathbb{R}^+\rightarrow \mathbb{C}\::\: f(\cdot) (\cdot)^{c-1}\in L^1(\mathbb{R}^+) \}$$
endowed with the norm
$$ \| f\|_{X_c} := \|f(\cdot) (\cdot)^{c-1} \|_1 = \int_0^{\infty} |f(u)|u^{c-1} du$$
(see \cite{BJ2}).
More generally, for $1<p<\infty$, we denote by $X^p_c$ the space of all  functions $f: \mathbb{R}^+\rightarrow \mathbb{C}$ such that 
$f(\cdot) (\cdot)^{c-1/p}\in L^p (\mathbb{R}^+)$; for $p=2$, see \cite{BJ4}. 
For $p=\infty$, we define $X^\infty_c$ as the space comprising 
all measurable functions $f : \mathbb{R}^+\rightarrow \mathbb{C}$ such that $\|f\|_{X^\infty_c}:= \sup_{x>0}x^{c}|f(x)| < \infty.$ 

For $f: \mathbb{R}^+ \rightarrow \mathbb{C},$ the Mellin translation operator $\tau_h^c$, $h \in \mathbb{R}^+$, is defined by
$$(\tau_h^c f)(x) \,:=\, h^c f(hx) \qquad(x\in \mathbb{R}^+).$$
Setting $\tau_h:= \tau^0_h,$ we have $(\tau_h^cf)(x) = h^c (\tau_hf)(x)$ and  $\|\tau_h^c f\|_{X_c} = \|f\|_{X_c}.$
\vskip0.3cm
The pointwise Mellin differential operator $\Theta_c$,
 $c \in \mathbb{R},$ is introduced as
$$
\Theta_cf(x) := x f'(x) + c f(x)\qquad(x \in \mathbb{R}^+)
$$
provided that $f'$ exists a.e.\ on $\mathbb{R}^+$ (see \cite{BJ1}). It is extended to order $r \in \mathbb{N}$ by defining recursively 
$\Theta^1_c := \Theta_c $, $\Theta^r_c := \Theta_c (\Theta_c^{r-1}).$
For convenience, we set $\Theta^r:= \Theta^r_0$ and $\Theta_c^0 := I$ with $I$ denoting the identity operator.
\vskip0,3cm
The Mellin transform of a function $f\in X_c$ is the linear and bounded operator defined by (see, e.g., \cite{MA}, \cite{GPS},  \cite{BJ2})
$$ M_c[f](s) \equiv [f]^{\wedge}_{M_c} (s) := \int_0^{\infty} u^{s-1} f(u) du \qquad(s=c+ it, t\in \mathbb{R}).$$
The inverse Mellin transform $M^{-1}_c[g]$ of a function $g \in L^1(c + i \mathbb{R})$ is given by
\begin{eqnarray*}
M^{-1}_c[g](x)  :=
 \frac{x^{-c}}{2 \pi}\int_{-\infty}^{\infty} g(c+it) x^{-it}dt \qquad(x \in \mathbb{R}^+),
\end{eqnarray*}
where  in general $L^p(c + i \mathbb{R}),$ for $p \geq 1,$ will mean the space of all functions 
$g:\, c+i \mathbb{R} \rightarrow \mathbb{C}$ such that 
$g(c +i\cdot) \in L^p(\mathbb{R}).$
\vskip0,3cm

The Mellin transform $M_c^2$ of $f \in X^2_c$ is given by (see \cite[Definition 2.1]{BJ4})
$$M_c^2[f](s) \equiv [f]^{\wedge}_{M_c^2} (s) = \mbox{l.i.m.}_{\rho \rightarrow +\infty}~\int_{1/\rho}^\rho f(u) u^{s-1}du
\qquad (s=c+it)$$
in the sense that
$$\lim_{\rho \rightarrow +\infty}\bigg\|M_c^2[f](c+it) - \int_{1/\rho}^\rho f(u) u^{c+it-1}du\bigg\|_{L^{2}(c + i \mathbb{R})} = 0.$$

Analogously, we define the inverse Mellin transform of a function $g \in L^2(c + i\mathbb{R})$ by 
$$\lim_{\rho \rightarrow +\infty}\bigg\|M^{2,-1}_c[g](x) - \frac{1}{2\pi}\int_{-\rho}^\rho g(c+it) x^{-c-it}dt\bigg\|_{X^2_c} = 0.$$
For any $f \in X^2_c$, there holds the inversion theorem
$$M^{2, -1}_c[M_c^2[f]](x) = f(x)\quad \quad (\hbox{a.e.\ on } \mathbb{R}^+);$$
see \cite[Theorem 2.9]{BJ4}.

For $f\in X_c \cap X^2_c,$ we have the important ``consistency'' property of the Mellin transform, namely $M_c[f](c+it) = M_c^2[f](c+it)$ for almost 
all $t \in \mathbb{R}$; see \cite[Lemma 2.5]{BJ4}.
\smallskip

In the following, we find it convenient to write $X^1_c:=X_c$, $M^1_c:= M_c$ and $M^{1,-1}_c:=M^{-1}_c$.
For $c \in \mathbb{R}$, $\sigma >0$ and $p\in \{1, 2\}$, the Mellin--Paley--Wiener spaces $B^p_{c,\sigma}$  
comprises all functions $f\in X^p_c \cap C(\mathbb{R}^+)$ such that
$[f]^\wedge_{M^p_c}(c+it) = 0$ for all $|t|>\sigma$ when $p=1$ and for almost all $|t|>\sigma$ when $p=2$.
The Mellin--Paley--Wiener space $B^p_{c, \sigma}$ is a subspace of the so-called Mellin inversion class 
$\mathcal{M}^p_c$ comprising all functions $f \in X^p_c \cap C(\mathbb{R}^+)$ such that $[f]^\wedge_{M^p_c} \in L^1(c + i\mathbb{R})$;
see \cite{BBMS2}.


One of the fundamental results for our theory is the following Mellin--Poisson summation formula, which was established in \cite{BJ2}. 
Denoting by 
$$\mathscr{M}_c[f](c+ik) := \int_{e^{-\pi}}^{e^\pi}f(u) u^{c+ik-1}du \qquad (k \in \mathbb{Z})$$
the Mellin-Fourier coefficients of $f$, we have (see \cite[Theorem 7.3]{BJ2}): 
\begingroup
\setcounter{tmp}{\value{Theorem}}
\setcounter{Theorem}{0} 
\renewcommand\theTheorem{\Alph{Theorem}}
\begin{Theorem}\label{mpst}
Let $f \in X_c,$ where $c \in \mathbb{R},$ be continuous on $\mathbb{R}^+$ and let
$$\sum_{k \in \mathbb{Z}}|\mathscr{M}_c[f](c+ik)| < +\infty.$$
If the series
$$f^c(x):= \sum_{k\in \mathbb{Z}} f(e^{2\pi k})e^{2\pi kc} \qquad (x \in \mathbb{R}^+)$$
converges uniformly for $x \in [e^{-\pi}, e^{\pi}],$ then
\begin{eqnarray}\label{mpsf}
f^c(x) = \frac{1}{2\pi}\sum_{k \in \mathbb{Z}}\mathscr{M}_c[f](c+ik) x^{-ik-c} \qquad (x \in \mathbb{R}^+).
\end{eqnarray}
\end{Theorem}
\endgroup
\setcounter{Theorem}{\thetmp}
\vskip0,3cm
As a basic space in our considerations we recall the Mellin--Sobolev space $W^{r,p}_c(\mathbb{R}^+)$. For
$r \in \mathbb{N}$ and $p\geq 1,$ it is defined as the set of all functions 
$f \in X^p_c$ for which there exists a function $g \in C^{(r-1)}(\mathbb{R}^+)$ such that $f=g$ a.e.\ and $g^{(r-1)}\in AC_{\rm{loc}}(\mathbb{R}^+)$ with   
$ \Theta^r_cg\in X^p_c.$
Here $AC_{\rm{loc}}(\mathbb{R}^+)$ denotes the space of all locally absolutely continuous functions on $\mathbb{R}^+.$
Briefly we may write
$$ W_c^{r,p}(\mathbb{R}^+) \,=\, \{ f\in X^p_c : \Theta^r_c f \in X^p_c \}.$$

 As a consequence of Theorem \ref{mpst}, we obtain the following statement, which is substantial for our purposes (see \cite[Corollary 7.5]{BJ2}).
\begin{cor}\label{mpst2}
 Let $f\in W_c^{1,1}(\mathbb{R}^+) $ be a continuous function on $\mathbb{R}^{+}$ such that for some $\sigma>0$
 $$ \sum_{k\in \mathbb{Z}} |[f]_{M_c}^{\wedge} (c+ 2\pi i \sigma k)| < \infty. $$
 Then 
 \begin{eqnarray}\label{mpsf2}
 \sum_{k\in \mathbb{Z}} f(e^{k/\sigma}) e^{kc/\sigma} = \sigma \sum_{k\in \mathbb{Z}} [f]_{M_c}^{\wedge}(c+ 2\pi i \sigma k). 
 \end{eqnarray}
 \end{cor}
 \vskip0,3cm
 If $f\in B^1_{c,2\pi\sigma}$, then the above equality reduces to 
 $$ \sigma  [f]_{M_c}^{\wedge} (c) = \sum_{k\in \mathbb{Z}} f(e^{k/\sigma}) e^{kc/\sigma} $$
 or, equivalently,
 \begin{eqnarray}\label{qf}
 \int_0^{\infty}  f(x) x^{c-1} dx = \frac{1}{\sigma}\sum_{k=-\infty}^{\infty} f(e^{k/\sigma}) e^{kc/\sigma},
 \end{eqnarray} 
 which can be interpreted as an exact quadrature formula. 

 If  $f \not \in B^1_{c,2\pi\sigma}$, this quadrature formula is no longer exact but it holds with a remainder term. Indeed we can rewrite 
(\ref{mpsf2}) as
 $$  \int_0^\infty f(x) x^{c-1} dx - \frac{1}{\sigma}\sum_{k=-\infty}^{\infty} f(e^{k/\sigma}) e^{kc/\sigma} = 
 - \sum_{k\in \mathbb{Z}, k\not =0} [f]_{M_c}^{\wedge} (c+ 2\pi i \sigma k).$$
 Thus, defining
\begin{equation}\label{def_remainder}
 R_{c,\sigma} [f] := -\sum_{k\in \mathbb{Z}, k\not =0} [f]_{M_c}^{\wedge} (c+ 2\pi i \sigma k),
\end{equation}
 we obtain the approximate quadrature formula 
 \begin{eqnarray}\label{ampsf}
\int_0^\infty f(x) x^{c-1} dx \,=\, \frac{1}{\sigma} \sum_{k=-\infty}^{\infty} f(e^{k/\sigma}) e^{kc/\sigma} +  R_{c,\sigma} [f].
  \end{eqnarray}

 In the next section we will study estimates for the remainder $R_{c,\sigma} [f]$ when the function $f$ belongs to suitable function spaces.


\subsection{Polar analytic functions, Mellin--Hardy spaces and Paley--Wiener type results}
\label{Sec:2.2}

Let $\mathbb{H}:= \{(r,\theta) \in \mathbb{R}^+ \times \mathbb{R}\}$ be the right half-plane and let ${\cal D}$ be a domain in $\mathbb{H}.$ 
For $a>0$ denote by $\mathbb{H}_a$  the set
$$\mathbb{H}_a:=\{(r,\theta) \in \mathbb{R}^+ \times ]-a,a[\}$$
and call it a strip in $\mathbb{H}$.

We begin with the following notion first introduced in \cite[Definition 3.1]{BBMS3}.
%
\begin{Definition}\label{def1}
We say that $f:{\cal D}\rightarrow \mathbb{C}$ is {\rm polar-analytic} on ${\cal D}$ if for any $(r_0, \theta_0) \in {\cal D}$ the limit
$$\lim_{(r,\theta) \rightarrow (r_0, \theta_0)}\frac{f(r, \theta) - f(r_0, \theta_0)}{re^{i\theta} - r_0e^{i\theta_0}} =: (D_{{\rm pol}}f)(r_0, \theta_0)$$
exists and is the same howsoever $(r, \theta)$ approaches $(r_0, \theta_0)$ within ${\cal D}.$
\vskip0,3cm
For a polar-analytic function $f$, we define the polar Mellin derivative as 
\begin{equation}\label{def1eq1}
\Theta_cf(r,\theta) := re^{i \theta}(D_{{\rm pol}}f)(r, \theta) + cf(r,\theta).
\end{equation}
\end{Definition}
\vskip0,3cm
\begin{Remark}
{\rm This modified notion of analyticity, arising by treating polar coordinates as 
cartesian coordinates,  leads naturally to  the classical Cauchy-Riemann
equations when written in their  polar form, i.e.
 $f = u + iv$ with $u,v: {\cal D}\rightarrow \mathbb{R}$ is polar-analytic on ${\cal D}$ if and only if $u$ and $v$ have continuous partial 
derivatives on ${\cal D}$ that satisfy the differential equations}
\begin{eqnarray}\label{CRE}
\frac{\partial u}{\partial \theta} = - r \frac{\partial v}{\partial r}\,,\qquad 
\frac{\partial v}{\partial \theta} = r \frac{\partial u}{\partial r}\,.
\end{eqnarray}
\end{Remark}
 For the derivative $D_{{\rm pol}}$, we easily find that 
$$(D_{{\rm pol}}f)(r, \theta) = e^{-i\theta}\bigg[\frac{\partial}{\partial r}u(r, \theta) + i \frac{\partial}{\partial r}v (r, \theta) \bigg] = 
\frac{e^{-i\theta}}{r}\bigg[\frac{\partial}{\partial \theta}v (r, \theta) - i \frac{\partial}{\partial \theta}u (r, \theta) \bigg].$$
Also note that $D_{{\rm pol}}$ is the ordinary differentiation on $\mathbb{R}^+.$ More precisely, if $\varphi (\cdot) := f(\cdot, 0)$, 
then $(D_{{\rm pol}}f)(r,0) = \varphi'(r).$
Moreover, for $\theta = 0$ we recover from (\ref{def1eq1}) the known formula  
$$\Theta_c\varphi(r) = r\varphi'(r) + c\varphi (r).$$

When $g$ is an entire function, then $f: (r, \theta) \mapsto g(re^{i\theta})$ defines a function $f$ on $\mathbb{H}$ that is polar-analytic 
and $2\pi$-periodic with respect to $\theta.$ The converse is also true. However, there exist polar-analytic functions on $\mathbb{H}$ 
that are not $2\pi$-periodic with respect to $\theta.$ A simple example is the function $L(r, \theta):= \log r + i\theta,$ which is easily seen 
to satisfy the differential equations (\ref{CRE}).
\vskip0,3cm
\begin{Definition}\label{MBS}
For $c \in \mathbb{R},$ $T>0$ and $p \in [1, +\infty[$ the {\rm Mellin--Bernstein space} $\mathscr{B}^p_{c,T}$ comprises all functions 
$f: \mathbb{H}\rightarrow \mathbb{C}$ with the following properties:
\begin{enumerate}
\item[(i)] $f$ is polar-analytic on $\mathbb{H};$
\item[(ii)] $f(\cdot, 0) \in X^p_{c};$
\item[(iii)] there exists a positive constant $C_f$ such that
$$|f(r,\theta)| \leq C_fr^{-c}e^{T|\theta|}\qquad ((r, \theta) \in \mathbb{H}).$$
\end{enumerate}
\end{Definition}
\vskip0,3cm
The Paley--Wiener theorem for the space $\mathscr{B}^2_{c,T}$ is the following statement (see \cite[Theorem 4]{BBMS3}):
\begingroup
\setcounter{tmp}{\value{Theorem}}
\setcounter{Theorem}{1} 
\renewcommand\theTheorem{\Alph{Theorem}}
\begin{Theorem}[Paley-Wiener]\label{pw}
A function $\varphi \in X^2_c$ belongs to the Mellin--Paley--Wiener space $B^2_{c,T}$ if and only if there exists a function $f\in \mathscr{B}^2_{c,T}$ such that 
$f(\cdot, 0) = \varphi(\cdot).$
\end{Theorem}
\endgroup
\setcounter{Theorem}{\thetmp}
\vskip0,4cm
For $c \in \mathbb{R}$ and  
$p \in [1,+\infty[$,  we recall (\cite{BJ1}, \cite{BJ4}) that the norm 
in $X^p_c$ is defined by
$$\|\varphi\|_{X^p_c}= \bigg(\int_0^{\infty} |\varphi(r)|^p r^{cp-1}dr\bigg)^{1/p}.$$
The Mellin--Hardy spaces for polar-analytic functions are defined as follows (see \cite{BBMS3} for the definition and properties).
\begin{Definition} \label{Hardy}
Let $a, c, p \in \mathbb{R}$ with $a>0$ and $p\geq 1.$ The {\rm Mellin--Hardy space} $H^p_c(\mathbb{H}_a)$ comprises all functions $f : \mathbb{H}_a \rightarrow \mathbb{C}$ that satisfy the following 
conditions:
\begin{enumerate}
\item[(i)] $f$ is polar-analytic on $\mathbb{H}_a;$
\item[(ii)] $f(\cdot, \theta) \in X^p_c$ for each $\theta \in {]}-a,a{[};$
\item[(iii)] there holds
$$\|f\|_{H^p_c(\mathbb{H}_a)} := \sup_{0<\theta<a}\bigg(\frac{\|f(\cdot, \theta)\|^p_{X^p_c} + \|f(\cdot, -\theta)\|^p_{X^p_c}}{2}\bigg)^{1/p} < +\infty.$$
\end{enumerate}
\end{Definition}
The links with the classical Hardy type spaces considered in \cite{BSS0} are explained in \cite{BBMS3}.

In order to state a generalization of the Mellin--Paley--Wiener theorem that characterizes the space of all functions such that their Mellin transform 
decays exponentially at infinity, we introduced in \cite[Defintion 3.4]{BBMS4} the following space $H^\ast_c(\mathbb{H}_a),$ which lies between two Mellin--Hardy spaces.
\vskip0,3cm
\begin{Definition}\label{intermediate}
Let $a,c \in \mathbb{R}$ with $a>0$ be fixed numbers. The class $H^\ast_c(\mathbb{H}_a),$ comprises all functions $f:\mathbb{H}_a \rightarrow \mathbb{C}$ with the following properties:
\begin{enumerate}
\item[(a)] $f$ is polar-analytic on $\mathbb{H}_a$;
\item[(b)] $\lim_{r \rightarrow 0^+} r^cf(r,0) = \lim_{r \rightarrow +\infty} r^cf(r,0) = 0$;
\item[(c)] for every $\varepsilon \in ]0,a[$ there exists a constant $K(f,\varepsilon)$ such that
$$|f(r, \theta)| \leq r^{-c}K(f,\varepsilon)$$
for all $(r,\theta) \in \mathbb{H}_{a-\varepsilon};$
\item[(d)] for every $\theta \in ]-a,a[$ and all $t \in \mathbb{R},$
$$I_c(f,\theta,t):= \lim_{R \rightarrow +\infty}\int_{1/R}^R f(r,\theta)r^{c+it-1}dr$$
exists and $|I_c(f,\theta, t)| \leq K(f)$ with a constant depending on $f$ only.
\end{enumerate}
\end{Definition}

We showed that  $ H_c^\ast(\mathbb{H}_a)$ becomes a normed linear space by endowing it with 
$$ \|f\|_{H_c^\ast(\mathbb{H}_a)}:= \sup \{ |I_c(f,\theta, t) | : \theta\in ]-a, a [, t\in \mathbb{R}\}.$$
As proved in \cite[Proposition 3]{BBMS4}, if  $f: \mathbb{H}_a \rightarrow \mathbb{C}$ satisfies conditions (a)--(c) of Definition \ref{intermediate}, then
$$ \lim_{r\rightarrow 0^+} r^c f(r, \theta)= \lim_{r\rightarrow +\infty} r^c f(r, \theta)=0$$
uniformly with respect to $\theta$ on all compact subintervals of $ ]-a, a [.$
\begin{Example}\label{ex1}
{\rm For $\mu >0,$ let $f(r):=\exp(-r^{\mu}).$ It is easily seen that
$$F(r,\theta):= \exp\bigg(-r^\mu e^{i\mu \theta}\bigg)$$
is a polar-analytic extension of $f$ to $\mathbb{H}.$ Furthermore it can be verified that $F \in H^\ast_c(\mathbb{H}_a)$ for all $c >0$ and $a \in ]0, \pi/(2\mu)[.$ However $F \not \in H^\ast_c(\mathbb{H}_a)$ when $a \geq \pi/(2\mu).$
}
\end{Example}
\vskip0,3cm

Subsequently, for a function $f \in H_c^\ast(\mathbb{H}_a),$ we set $\phi(r):= f(r,0)$ and 
$$M^\ast_c[\phi](c+it):= I_c(f(r,\cdot), 0,t) \qquad  (t \in \mathbb{R}).$$
The following generalization of the Paley--Wiener theorem for the Mellin transform holds (see \cite[Theorem 3.2]{BBMS4}):
\begingroup
\setcounter{tmp}{\value{Theorem}}
\setcounter{Theorem}{2} 
\renewcommand\theTheorem{\Alph{Theorem}}
\begin{Theorem} \label{expopw}
A continuous function $\phi:\mathbb{R}^+ \rightarrow \mathbb{C}$ is the restriction to $\mathbb{R}^+$ of a function $f \in H^\ast_c(\mathbb{H}_a)$ if and only if $M^\ast_c[\phi]$ exists and 
\begin{eqnarray}\label{16}
|M^\ast_c[\phi](c+it)| \leq Ce^{-a|t|}\qquad (t \in \mathbb{R})
\end{eqnarray}
with a constant $C$ that may be taken to be $\|f\|_{H^\ast_c(\mathbb{H}_a)}.$
\end{Theorem}
\endgroup
\setcounter{Theorem}{\thetmp}

\section{Quadrature over the positive real axis}
\label{Sec:3}

In this section our chief aim is to estimate the series (\ref{def_remainder}) defining the remainder $R_{c,\sigma},$ a core part of our paper.
We first introduce a basic class of functions for which the Mellin--Poisson summation formula holds (see \cite[Definition 3.1]{SCH})
\begin{Definition}\label{MP}
Let $c \in \mathbb{R}$ and $s>0.$ Then the class $\mathcal{K}_c(s)$ comprises all functions $f \in X_c\cap C(\mathbb{R}^+)$ such that 
$$
\sum_{k=-\infty}^\infty |f(e^{k/\sigma})| e^{kc/\sigma} < +\infty \quad \hbox{ and } \quad
\sum_{k=-\infty}^\infty |[f]^\wedge_{M_c}(c+ i 2\pi k \sigma)| < +\infty
$$
for all $\sigma \geq s.$
\end{Definition}
The new results to be deduced here give estimates for the remainder of the quadrature formula (\ref{ampsf}) when
 the involved functions belong to the class $\mathcal{K}_c(s).$ 

We divide this section in two parts. In the first one we give estimates of the remainder in Mellin--Sobolev type spaces of integer order
and in the second one we discuss a fractional case.

\subsection{Estimates in Mellin--Sobolev spaces of integer order}
\label{Sec:3.1}

In the Mellin setting, the notions of modulus of smoothness and Lipschitz classes are introduced as in the classical case
(see, e.g., \cite[Chap.~2, Sects.~7 and 9]{DL}) except that the role of the difference operator $\Delta_h$ is now taken up by the
Mellin translation operator. More precisely, we
define the difference of integer order $r\in \mathbb{N}$ of a function $f \in X_c$  as
 $$ (\Delta_h^{r,c} f) (u):=  \sum_{j=0}^r (-1)^{r-j} \left(\begin{array}{c} r\\ j\end{array} \right) f(h^j u) h^{jc}.$$
It is needed for introducing 
 the modulus of smoothness 
$$ \omega_r(f, \delta, X_c) := \sup_{|\log h|\leq \delta} \|\Delta_h^{r,c} f\|_{X_c}, \quad \delta>0,$$
which in turn is used for defining the Mellin--Lipschitz classes
 $$\mbox{ Lip}_r (\alpha, X_c) := \{ f\in X_c\::\: \omega_r(f, \delta, X_c) = {\mathcal O}(\delta^\alpha), \, \delta \rightarrow 0^+ \},$$
where $\alpha \in ]0, r]$.
In this setting the remainder $R_{c,\sigma}[f]$ of the quadrature formula (\ref{ampsf})
can be estimated asymptotically as follows.
\begin{Theorem}\label{lip}
 Let $f\in \mbox{\rm Lip}_r (\alpha, X_c) \cap C(\mathbb{R}^+)$ with $r\geq 2$ being an integer and $1<  \alpha \leq r$. Then
\begin{eqnarray}\label{est1}
|R_{c,\sigma}[f]| = \mathcal{O}(\sigma^{-\alpha})\quad \quad (\sigma \rightarrow +\infty).
\end{eqnarray}
\end{Theorem}
{\bf Proof}. As shown in the proof of  \cite[Theorem 3]{BBMS2}, the hypotheses imply that
 $$ |[f]_{M_c}^{\wedge} (c + iv) | \,\leq\, \frac{1}{2^r}\,\omega_r\left(f, \frac{\pi}{|v|}, X_c\right).$$
 Therefore, by (\ref{def_remainder}), 
 $$  |R_{c, \sigma} [f] | \,\leq\, \sum_{k\in \mathbb{Z}, k\not =0} |[f]_{M_c}^{\wedge} (c+ 2\pi i \sigma k)| \,\leq \,
\sum_{k\in \mathbb{Z}, k\not =0}\frac{1}{2^r}\, \omega_r\left(f, \frac{1}{2\sigma |k|}, X_c\right).$$
Since $f \in \mbox{Lip}_r(\alpha, X_c)$,  there exists a constant $C>0$ such that for sufficiently large
 $\sigma>0$ we have
 \begin{align*}
 |R_{c, \sigma}[f| &\leq\, \frac{1}{2^{r-1} }\sum_{k =1}^\infty \omega_r\left(f, \frac{1}{2\sigma k}, X_c\right) \,\leq\,
 \frac{C}{2^{r-1}} \sum_{k=1}^\infty \frac{1}{(2\sigma k)^\alpha}\\[.7ex]
& =\, \frac{C}{2^{r+\alpha-1} \sigma^\alpha}\,\zeta(\alpha) \,= \,{\mathcal{O}}(\sigma^{-\alpha}) \qquad (\sigma\to +\infty),
 \end{align*} 
which implies (\ref{est1}). Here $\zeta(\cdot)$ denotes the Riemann zeta function. $\:\Box$
\vskip0,5cm
Under the existence of higher orders of derivatives, the previous asymptotic estimate can be improved.
\begin{Theorem}\label{sobolev}
Let $f \in W^{j,1}_c(\mathbb{R}^+)$ be such that $\Theta^j_cf \in \mbox{\rm Lip}_r(\alpha, X_c),$ where $\alpha\in {]}1, j{[}$. Then
\begin{eqnarray}\label{est2}
|R_{c,\sigma}[f]| = \mathcal{O}(\sigma^{-\alpha - j})\qquad (\sigma \rightarrow +\infty).
\end{eqnarray}
\end{Theorem}
{\bf Proof}. For $f \in W^{j,1}_c(\mathbb{R}^+)$ one has (see \cite[Proposition 6]{BJ1})
 $$[\Theta_c^j f]^{\wedge}_{M_c} (c+iv) \,=\, (-iv)^j [f]^{\wedge}_{M_c} (c+ iv)\qquad(v\in \mathbb{R})$$ 
or
 $$  [f]^{\wedge}_{M_c} (c+ iv) \,=\, (-iv)^{-j} [\Theta_c^j f]_{M_c}^{\wedge} (c+ iv), $$
 and so
 $$ |[f]^{\wedge}_{M_c} (c+ iv) | \,=\,  |v|^{-j }| [\Theta_c^j f]_{M_c}^{\wedge} (c+ iv) | \,\leq \, \frac{1}{2^r} |v|^{-j} 
\omega_r\left(\Theta^j_c f , \frac{\pi}{|v|}, X_c\right).$$
 Therefore
 \begin{align*}
 |R_{c,\sigma} [f]| &\leq\, \sum_{k\in \mathbb{Z}, k\not =0} |[f]_{M_c}^{\wedge} (c+ 2\pi i \sigma k)| \,\leq\, 
\frac{1}{2^r} \sum_{k\in \mathbb{Z}, k\not =0} \frac{1}{(2\pi \sigma k)^j}i\, \omega_r\left(\Theta^j_c f, \frac{\pi}{2 \pi \sigma |k|} X_c\right)\\[.7ex]
 &=\,\frac{1}{2^{r+j-1}(\pi\sigma)^j}\, \sum_{k=1}^\infty \frac{1}{k^j}\, \omega_r\left(\Theta^j_c f, \frac{1}{2  \sigma k}, X_c\right).
 \end{align*}
 Since $\Theta_c^j f \in \mbox{Lip}_r (\alpha, X_c)$, there exists a constant $C>0$ such that for sufficiently large $\sigma$ we have
\begin{align*}
 |R_{c,\sigma}[f]| &\leq\,  \frac{1}{2^{r+j-1}(\pi\sigma)^j} \,\sum_{k=1}^\infty \frac{C}{k^j(2\sigma k)^\alpha} \\[.7ex]
&=\,  \frac{C}{2^{r+j+\alpha-1}\pi^j \sigma^{j+\alpha}}\, \zeta(j+\alpha)\, = \,{\mathcal O}(\sigma^{-\alpha - j}) \qquad (\sigma \rightarrow +\infty),
 \end{align*}
which is (\ref{est2}). $\:\Box$


\subsection{Estimates in Mellin--Sobolev spaces of fractional order}
\label{Sec:3.2}

Now we discuss the fractional case. Let $\alpha >0$ be a fixed real number. In the following, powers of order $\alpha$  are defined with the 
help of the principal value of the logarithm. In particular, $(-1)^\alpha:= \exp(i\alpha\pi)$ and for $v\in\mathbb{R}$,
$$ (iv)^\alpha := |v|^\alpha \exp(i\alpha (\pi/2) \mbox{sgn} (v)) = i^\alpha v^\alpha.$$

For  $f:\mathbb{R}^+ \rightarrow \mathbb{C}$, the Mellin difference of order $\alpha$ is defined by the series
\begin{eqnarray*}
\Delta_{h}^{\alpha, c}f(x):= \sum_{j=0}^\infty \left(\begin{array}{l} \alpha\\j\end{array}\right) (-1)^{\alpha - j}\tau^c_{h^j}f(x)\qquad (x>0, \,h >0).
\end{eqnarray*}
It is known that for $f \in X_c$ and any $h>0$ the fractional difference exists a.e.\ and
$$\|\Delta_{h}^{\alpha, c}f\|_{X_c} \leq \|f\|_{X_c}\sum_{j=0}^\infty \left|\left(\begin{array}{l} \alpha\\j\end{array}\right)\right|$$
with a convergent series on the right-hand side
(see \cite[Proposition 1]{BBM2} and \cite[page 124]{BW}).

\begin{Definition}\label{fractional}
We define the (pointwise) Mellin fractional derivative of order $\alpha >0$ at a point $x$ by the limit 
$$\lim_{h\rightarrow 1}\frac{\Delta_{h}^{\alpha, c}f(x)}{(h-1)^\alpha}=: \Theta^\alpha_c f(x)$$
provided that it exists.
\end{Definition}
In \cite{BBM2} the strong fractional derivative for functions $f \in X_c$ is introduced 
as $\mbox{s-}\Theta^\alpha_cf:=g$ with $g$ satisfying the formula
$$\lim_{h \rightarrow 1}\left\|\frac{\Delta_h^{\alpha, c}f(\cdot)}{(h-1)^\alpha} - g(\cdot)\right\|_{X_c} = 0.$$ 
Moreover, 
$$[\mbox{s-}\Theta_c^\alpha f]^\wedge_{M_c}(c+iv) \,=\, (-iv)^\alpha[f]^\wedge_{M_c}(c+iv) \qquad (v \in \mathbb{R}).$$
Hence, if $f \in X_c$ has a strong fractional derivative of order $\alpha$ and also a pointwise fractional derivative $\Theta^\alpha_cf \in X_c,$ then 
\begin{eqnarray}\label{fractransf}
[\Theta^\alpha_cf]^\wedge_{M_c}(c+iv) = [\mbox{s-}\Theta_c^\alpha f]^\wedge_{M_c}(c+iv)=(-iv)^\alpha[f]^\wedge_{M_c}(c+iv).
\end{eqnarray}
For $p\in \{1, 2\}$, we now define the  Mellin--Sobolev space $W^{\alpha,p}_c(\mathbb{R}^+)$ of fractional order $\alpha$ as 
$$W^{\alpha,p}_c(\mathbb{R}^+) \,:=\,\left\{f \in X^p_c \::\: \Theta^\alpha_cf  ~\mbox{exists a.e. and }  ~ \Theta^\alpha_cf \in X^p_c \right\}.$$

We recall the following proposition (see \cite[Proposition~2.1]{BBMS4}):
\begin{prop} \label{fracrep} 
Let $p\in\{1, 2\}$, $c\in\mathbb{R}$, $f \in \mathcal{M}^p_c$, and let $\alpha >0$. 
Define $\varphi(v):=[f]^\wedge_{M^p_c}(c+iv)$.
If $v^\alpha \varphi (v) \in L^1(\mathbb{R})$, 
then $\Theta^\alpha_cf$ exists and 
$$\Theta^\alpha_cf(x) = \frac{(-i)^\alpha}{2 \pi}\int_{-\infty}^{\infty} v^\alpha \varphi (v) x^{-c-iv}dv \qquad (x>0).$$
\end{prop}
\vskip0,5cm
For the remainder $R_{c,\sigma}[f]$, we can now deduce the following asymptotic estimate.

\begin{Theorem}\label{fracsobolev}
Let $f \in \mathcal{M}^1_c \cap W^{\alpha, 1}_c(\mathbb{R}^+),$ where $\alpha >0$ and $c \in  \mathbb{R}.$ 
If $\Theta^\alpha_c f \in \mbox{\rm Lip}_r(\beta, X_c)$ for some $\beta >0$ such that $\alpha + \beta > 1$ and $r \in \mathbb{N},$ then 
$$|R_{c,\sigma}[f]| = \mathcal{O}(\sigma^{-\alpha -\beta})\qquad (\sigma \rightarrow +\infty).$$
\end{Theorem}
{\bf Proof}. By using (\ref{fractransf}), the proof can be performed analogously to that of Theorem~\ref{sobolev}. $\:\Box$
\medskip

For the space of all functions $f\,:\,\mathbb{R}^+ \to \mathbb{C}$ that have a representation
\begin{equation}\label{repr_Mellin} 
f(x)\,=\,\int_{-\infty}^\infty \varphi(v) x^{-c-iv} dv,
\end{equation}
where $\varphi\in L^1(\mathbb{R})\cap L^q(\mathbb{R})$, we introduced
 in \cite[page~301]{BBMS2} a metric $\dist_q$, induced by the $L^q(\mathbb{R})$ norm of $\varphi$. This metric is of special interest in
Mellin analysis. Indeed, if $f$ belongs to a Mellin inversion class, then a representation (\ref{repr_Mellin}) holds with $\varphi$ being a
Mellin transform of $f$.
 Here we want to employ this concept for estimates of the remainder $R_{c,\sigma}[f]$ in terms of the distance of
$f$ (or its Mellin derivatives) from the Mellin--Paley--Wiener space $B_{c,\sigma}^1$. We shall need this metric only for the Mellin inversion class
$\mathcal{M}_c^1$ and with $q=\infty$. In this case it can be introduced as
$$\dist_\infty(f,g)\,:=\, \sup_{v\in\mathbb{R}}\left|[f]_{M_c}^\wedge(c+iv) - [g]_{M_c}^\wedge(c+iv)\right|\qquad (f, g \in \mathcal{M}_c^1).$$
As a special case of \cite[Theorem~1, Corollary~1]{BBMS2}, we obtain the following formulae for distances from $B_{c,\sigma}^1$.
\begin{cor}\label{representation}
Let $f\in\mathcal{M}_c^1$. Then
$$ \dist_\infty(f, B_{c,\sigma}^1)\,=\, \sup_{|v|\ge\sigma}\left|[f]_{M_c}^\wedge(c+iv)\right|.$$
If, in addition, $\alpha>0$ and $(\cdot)^\alpha [f]_{M_c}^\wedge(c+i\cdot)\in L^1(\mathbb{R})\cap C(\mathbb{R})$, then $\Theta_c^\alpha f$ exists
and
$$ \dist_\infty(\Theta_c^\alpha f, B_{c,\sigma}^1)\,=\,  \sup_{|v|\ge\sigma}\left|v^\alpha[f]_{M_c}^\wedge(c+iv)\right|.$$
\end{cor}
\smallskip

Concerning the remainder in (\ref{ampsf}), we can now state the following result.
\begin{Theorem}\label{estdistance}
Let $f \in \mathcal{M}^1_c \cap W^{\alpha, 1}_c(\mathbb{R}^+)$ with $\alpha >1.$ Then 
\begin{equation}\label{Sob_bound}
|R_{c,\sigma}[f]| \,\leq\, \frac{2 \zeta(\alpha)}{(2\pi \sigma)^\alpha}\,\dist_\infty(\Theta_c^\alpha f, B^1_{c, 2\pi \sigma}).
\end{equation}
This estimate is best possible in the sense that it is no longer true if the right-hand side is multiplied by a positive constant less than $1$.
\end{Theorem}
{\bf Proof}.
Obviously, we have
\begin{eqnarray*}
|R_{c,\sigma}[f]| &\leq& \sum_{k \in \mathbb{Z}, k\neq 0}\frac{1}{|2\pi k\sigma|^\alpha}|(2\pi k\sigma)^\alpha [f]^\wedge_M(c+i2\pi k\sigma)|\\[.7ex]
&\leq& \frac{2\zeta(\alpha)}{(2\pi\sigma)^\alpha} \sup_{|v| \geq 2\pi\sigma} |v^\alpha [f]^\wedge_M(c+iv)|
\,\leq\, \frac{2\zeta(\alpha)}{(2\pi\sigma)^\alpha}\dist_\infty(\Theta^\alpha_cf, B^1_{c, 2\pi\sigma}),
\end{eqnarray*}
which is (\ref{Sob_bound}). Example~4 in Section~\ref{Sec:5} shows that the right-hand side cannot be diminished in the described sense.
 $\:\Box$  

\begin{Remark}
{\rm
For the estimate corresponding to (\ref{Sob_bound}) in the Fourier case, the authors \cite[Sec.~9]{BSS} have shown by a sophisticated construction
that for each $\sigma >0$ and each $\alpha >1$ there exists a function $\phi_{\alpha,\sigma}$ for which equality is attained. 
One can show that this stronger form of
sharpness also holds for (\ref{Sob_bound}) by modifying the extremal function $\phi_{\alpha,\sigma}$ suitably. 
}
\end{Remark}

The estimate of Theorem~\ref{estdistance} also holds when the remainder functional on the left-hand side is applied to a Mellin translation
of $f$. 
\begin{cor}\label{trans}
Under assumptions of Theorem \ref{estdistance} we have 
$$|R_{c,\sigma}[\tau^c_hf]| \,\leq\, \frac{2 \zeta(\alpha)}{(2\pi \sigma)^\alpha}\dist_\infty(\Theta_c^\alpha f, B^1_{c, 2\pi \sigma})$$
for every $h>0.$
\end{cor}
{\bf Proof}. By using the formula (see \cite[Lemma 3]{BJ1})
$$[\tau^c_hf]^\wedge_{M_c}(c+it) = h^{-it}[f]^\wedge_{M_c}(c+it),$$ 
we can procede as in the previous proof. $\:\Box$

\section{Equivalence theorems}
\label{Sec:4}

The results of Section 3 show that high regularity of the function $f$ yields rapid convergence to zero of the remainders of the quadrature 
formula (\ref {ampsf}) as $\sigma \rightarrow +\infty.$ One may ask if the converse is also true. Does rapid convergence to zero of the remainders 
imply high regularity of the function? The answer is {\it no} for a simple reason. Due to certain symmetries of the function $f,$ the remainders 
may be zero even if the graph of $f$ is very erratic as far as regularity is concerned. In  fact, we shall see subsequently that $f$ can be split 
into two parts, called the Mellin-even and the Mellin-odd part, such that the remainders of the Mellin-odd part are always zero. In order to 
eliminate this phenomenon, there are two possibilities: One may either study the remainders in terms of the Mellin-even part only or, alternatively, 
one may require that a certain rate of convergence not only holds for $f$ itself but also under the action of the Mellin translation operator 
$\tau_h^c$ on $f$ with values of $h$ near to $1$, which will destroy possible symmetries. Both approaches will be considered in this section. We mention 
that some basic results of this kind were obtained in \cite{DRS}, \cite{RS2} for quadrature formulae over the whole real axis and in \cite{SCH} for the 
positive real axis.

First we characterize the speed of convergence of the remainders in terms of distances from Mellin--Bernstein spaces. Then, using certain results 
proved in \cite{BBMS4}, which characterize distances from Mellin--Bernstein spaces by functions spaces guaranteeing a certain regularity,  
we can deduce a characterization of the speed of convergence by specific function spaces. 

%

\subsection{Characterization of speed of convergence by distances}
\label{Sec:4.1}
\subsubsection{Case of the Mellin-even part}
For $f:\mathbb{R}^+ \rightarrow \mathbb{C}$, we define (see \cite[Section IV]{SCH})
$$f_{c+}(x):= \frac{1}{2}(x^cf(x) + x^{-c}f(1/x)), \quad f_{c-}(x):=  \frac{1}{2}(x^cf(x) - x^{-c}f(1/x)) \quad (x>0),$$
so that $f(x) = x^{-c}(f_{c+}(x) + f_{c-}(x))$ for every $x>0.$
We call $f_{c+}$ and $f_{c-}$ the $c$-Mellin-even and the $c$-Mellin-odd part of $f,$ respectively. 
Moreover, for any function $f \in X_c$, we have $f_{c+}, f_{c-} \in X_0$ and, for $t \in \mathbb{R},$
$$[f_{c+}]^\wedge_{M_0}(\pm it) = \frac{1}{2}\bigl([f]^\wedge_{M_c}(c+it) + [f]^\wedge_{M_c}(c-it) \bigr)$$
and
$$[f_{c-}]^\wedge_{M_0}(it) = \frac{1}{2}\bigl([f]^\wedge_{M_c}(c+it) - [f]^\wedge_{M_c}(c-it) \bigr).$$
As was observed in \cite[Proposition 4.1]{SCH}, the remainders in (\ref{ampsf}) satisfy the following relations.
\begin{prop}\label{evenodd}
Let $f \in \mathcal{K}_c(s).$ Then
$$R_{c,\sigma}[f] = R_{0,\sigma}[f_{c+}]$$
and $R_{0,\sigma}[f_{c-}] = 0$ for all $\sigma \geq s.$
\end{prop}
This leads us to the representation
\begin{eqnarray}\label{representation3}
R_{c,\sigma}[f] = - 2\sum_{k=1}^\infty [f_{c+}]^\wedge_{M_0}(i2\pi k \sigma).
\end{eqnarray}
Using the so-called M\"{o}bius function $\mu: \mathbb{N}\rightarrow \{-1,0,1\}$, defined by
\begin{eqnarray*}
\mu(k) = \left\{\begin{array}{lll} 1, \quad & k=1,\\ (-1)^n, \quad & k=p_1\cdots p_n ,~\mbox{with distinct primes}~p_i,\\
0, \quad & k,~\mbox{is divisible by a square of a prime,}
\end{array}\right.
\end{eqnarray*}
we can invert formula (\ref{representation3}). In fact, we have (see \cite[Lemma 5.1]{SCH}):
\begin{lemma}\label{inverseformula}
Let $f \in \mathcal{K}_c(s).$  If 
\begin{equation}\label{condition}
\sum_{n=1}^\infty |R_{c,n\sigma}[f]| < +\infty
\end{equation}
for some $\sigma \geq s,$ then
\begin{eqnarray}\label{inverse}
[f_{c+}]^\wedge_{M_0}(i2\pi n \sigma)  = -\frac{1}{2}\sum_{k=1}^\infty \mu(k) R_{c,nk\sigma}[f]
\end{eqnarray}
for all $n \in \mathbb{N}.$ In particular
\begin{eqnarray}\label{inverse2}
|[f_{c+}]^\wedge_{M_0}(i2\pi \sigma) | \leq  \frac{1}{2}\sum_{k=1}^\infty  |R_{c,k\sigma}[f]|.
\end{eqnarray}
\end{lemma}
\vskip0,3cm
Next we describe an important situation where condition (\ref{condition}) is satisfied.  
\begin{lemma}\label{lemma3}
Let $f \in \mathcal{K}_c(s)$  and $\alpha >1.$ Under the assumptions of Theorem \ref{estdistance} we have
$$\sum_{n=1}^\infty |R_{c,n\sigma}[f]| < +\infty.$$
\end{lemma}
{\bf Proof}. Taking into account that for every $n \geq 1$ one has $B^1_{c,2 \pi \sigma}  \subset B^1_{c,2 \pi n \sigma},$ 
we deduce from the  estimate in Theorem \ref{estdistance} that
$$||R_{c,n\sigma}[f]| \leq \frac{2 \zeta(\alpha)}{(2\pi n \sigma)^\alpha} \dist_\infty(\Theta^\alpha_cf, B^1_{c, 2\pi \sigma}).$$
Now the assertion follows immediately. $\:\Box$
\vskip0,4cm
We are now ready for the main result of this subsection.

\begin{Theorem}\label{equiv2}
Let $f \in \mathcal{K}_c(s).$ Let $\lambda$ be a non-negative, nonincreasing function 
on an interval $[t_0, +\infty[$ with $t_0>0$, and let the assumptions of Theorem~\ref{estdistance} be satisfied.
Then, for $\alpha >1$, the following statements are equivalent:
\begin{itemize}
\item[(i)] $\;\dist_\infty(\Theta^\alpha f_{c+}, B^1_{0, 2\pi \sigma}) = \mathcal{O}(\lambda(\sigma))\quad (\sigma \rightarrow +\infty)$;\\[-1ex]

\item[(ii)] $\;R_{c,\sigma}[f] = \mathcal{O}(\sigma^{-\alpha}\lambda(\sigma)) \quad (\sigma \rightarrow +\infty).$
\end{itemize}
\end{Theorem}
{\bf Proof}. Suppose that (i) holds. Since $f_{c+}$ satisfies the same assumptions as $f,$ but with $c=0,$ 
we see that (ii) follows from Theorem \ref{estdistance} in conjunction with Proposition~\ref{evenodd}.

Conversely, suppose that (ii) holds. Then there exist a constant $C>0$ and $\sigma_0 > t_0$ such that for every $\sigma \geq \sigma_0$ we have, 
by Proposition \ref{evenodd},
$$|R_{c,\sigma}[f_{c+}]| \,= \,|R_{c,\sigma}[f]| \,\leq\, C\,\frac{\lambda(\sigma)}{\sigma^\alpha} \qquad (\sigma \geq \sigma_0).$$
Since Lemma~\ref{inverseformula} is applicable, we obtain with the help of (\ref{inverse2}) that
$$|[f_{c+}]^\wedge_{M_0}(i2\pi \sigma)| \,\leq \,\frac{1}{2}\sum_{k=1}^\infty |R_{c,k\sigma}[f]| \,\leq \,
\frac{C}{2}\sum_{k=1}^\infty \frac{\lambda(k\sigma)}{(k\sigma)^\alpha}.$$
Therefore,  by the hypotheses on the function $\lambda,$
$$|(2\pi\sigma)^\alpha [f_{c+}]^\wedge_{M_0}(i2\pi \sigma)| \,\leq \,\frac{C}{2}\,(2\pi)^\alpha\zeta(\alpha) \lambda(\sigma),$$
for every $\sigma \geq \sigma_0.$ This implies that
$$\sup_{|v| \geq 2\pi \sigma}\{|v|^\alpha |[f_{c+}]^\wedge_{M_0}(iv)|\} \,\leq \,\frac{C}{2}\,(2\pi)^\alpha\zeta(\alpha)\lambda(\sigma) 
\qquad (\sigma \geq \sigma_0),$$
and (i) follows by Corollary~\ref{representation}.  $\:\Box$ 
\vskip0,3cm
\subsubsection{Case of the translated function $\tau^c_hf$}\label{subsub4.1.2}
Now we establish a characterization of the speed of convergence  of the remainder for the translated function $\tau_h^cf$  in terms of the distance 
$\dist_\infty(f, B^1_{c, 2\pi \sigma}).$
As a counterpart of Lemma~\ref{inverseformula} for remainders of translated functions, we state the following auxiliary result.
 
\begin{lemma}\label{lem_shift}
Let $f\in \mathcal{K}_c(s)$. Then
$$
\left|[f]^\wedge_{M_c}(c\pm i2\pi\sigma)\right|\,\le\, 
\sup\left\{\left|R_{c,\sigma}\left[\tau_h^c f\right]\right| \::\: e^{-1/(2\sigma)}\le h\le e^{1/(2\sigma)}\right\}
$$
for $\sigma>s.$
\end{lemma}
{\bf Proof.}\, 
For short, we write the Mellin translation as $f_h:=\tau^c_hf.$ We recall that
$$
[f_h]^\wedge_{M_c}(c+it) = h^{-it}[f]^\wedge_{M_c}(c+it), 
$$
which implies by (\ref{def_remainder}) that
$$R_{c,\sigma}[f_h] = -\sum_{k \in \mathbb{Z}, k \neq 0} h^{-i2\pi k \sigma}[f]^\wedge_{M_c}(c + i2\pi k\sigma).$$
Setting $x:= h^{2\pi \sigma},$ we may rewrite this equation as
$$-\frac{x^{-c}}{2\pi} R_{c,\sigma}[f_{x^{1/(2\pi \sigma)}}] = \frac{x^{-c}}{2\pi}\sum_{k \in \mathbb{Z}, k \neq 0}  [f]^\wedge_{M_c}(c + i2\pi k\sigma)x^{-ik}.$$
The right-hand side is a Mellin--Fourier series in $x$
(see \cite[Section~5]{BJ2}).  When $h$ traverses the interval $[e^{-1/(2\sigma)}, e^{1/(2\sigma)}],$ then $x$ traverses the interval $[e^{-\pi}, e^{\pi}].$
Therefore, by \cite[Theorem~3.2]{BJ4},
we have
\begin{eqnarray}\label{18}
\sum_{k \in \mathbb{Z}, k \neq 0} \big|[f]^\wedge_{M_c}(c + i2\pi k\sigma)\big|^2 \,=\,
\frac{1}{2\pi} \int_{e^{-\pi}}^{e^\pi} |R_{c,\sigma}[f_{x^{1/(2\pi \sigma)}}]|^2\,\frac{dx}{x}.
\end{eqnarray}
From this we deduce that
$$ \sum_{k \in \mathbb{Z}, k \neq 0} \big|[f]^\wedge_{M_c}(c + i2\pi k\sigma)\big|^2 \,\le\,
\sup\left\{\left|R_{c,\sigma}\left[f_h\right]\right|^2 \::\: e^{-1/(2\sigma)}\le h\le e^{1/(2\sigma)}\right\}.
$$
Now, considering on the left-hand side the terms for $k= \pm 1$ only, we obtain the assertion. $\:\Box$ 

As an analogue of Theorem~\ref{equiv2}, we have:
\begin{Theorem}\label{equiv1}
Let $f \in \mathcal{K}(s).$  
Let $\lambda$ be a non-negative, non-increasing function on $[t_0, +\infty[$ with $t_0>0$, and let the assumptions of Theorem~\ref{estdistance} be satisfied.
Then, for $\alpha >1$, the following assertions are equivalent:
\begin{enumerate}
\item[(i)] $\;\dist_\infty(\Theta^\alpha_cf, B^1_{c, 2\pi \sigma}) = \mathcal{O}(\lambda(\sigma))$ as $\sigma \rightarrow +\infty;$\\[-1ex]

\item[(ii)] $\;R_{c,\sigma}[\tau^c_hf] = \mathcal{O}(\sigma^{-\alpha}\lambda(\sigma))$ as $\sigma \rightarrow +\infty$ uniformly for $h \in [e^{-1/(2\sigma)}, e^{1/(2\sigma)}].$
\end{enumerate}
\end{Theorem}
{\bf Proof}. If (i) holds, then  Corollary \ref{trans} implies immediately that (ii) is true.

Conversely, if (ii) holds, then, by Lemma~\ref{lem_shift},
$$[f]^\wedge_{M_c}(c \pm  i 2\pi\sigma)| \,= \,\mathcal{O}\left(\sigma^{-\alpha}\lambda(\sigma)\right)$$
as $\sigma\to +\infty$, which implies that
$$\dist_\infty(\Theta^\alpha_cf, B^1_ {c, 2\pi \sigma})  
\,=\, \sup_{|v| \geq 2\pi \sigma}|v^\alpha [f]^\wedge_{M_c}(c + iv)| \,= \, \mathcal{O}(\lambda(\sigma))\qquad (\sigma\to +\infty),$$
and so (i) holds. $\:\Box$

\subsection{Characterizations of distances by function spaces}
\label{Sec:4.2}

Our aim is to determine function spaces which guarantee a prescribed rate of convergence of remainders.
For this we will employ results obtained in \cite{BBMS4} which characterize distances from Mellin--Paley--Wiener spaces by specific function spaces.
For the reader's convenience, we reproduce the former theorems as propositions. We start with \cite[Theorem~4.1]{BBMS4}.

\begin{prop}\label{exactness}(Exactness)
Let $f \in \mathcal{M}^1_c.$ Then
$$f \in B^1_{c, \sigma}\: \Longleftrightarrow \: \dist_\infty(f, B^1_{c,\rho} )= 0\quad (\forall\: \rho>\sigma).$$
\end{prop}

The next result can be deduced from Theorem~\ref{expopw}; for details see \cite[Theorem~4.2]{BBMS4}. We reproduce it in a concise form. 
\begin{prop}\label{exprate}(Exponential rate) 
Let $\varphi \in \mathcal{M}^1_c.$ Then
$$
\left.
\begin{array}{c}
f \in H^\ast_c(\mathbb{H}_a)\\[.8ex]
f(\cdot,0)\equiv \varphi
\end{array} 
\right\}
\: \Longleftrightarrow \: \dist_\infty(\varphi, B^1_{c,\sigma}) = \mathcal{O}(e^{-a \sigma})\quad  (\sigma \rightarrow +\infty).$$
\end{prop}
\vskip0,2cm

In \cite{BBMS2} we proved that functions from  a Mellin--Sobolev space of order $r$ have distances from $B^1_{c,\sigma}$ 
that behave like $\mathcal{O}(\sigma^{-r}).$ However, as remarked in \cite{BBMS4}, the familiar Mellin--Sobolev spaces are not appropriate 
for characterizing this rate of convergence. The following modified space will accomplish the desired equivalence trivially.%
\footnote{The rate characterized by the Mellin--Sobolev space $W_0^{r,2}(\mathbb{R}^+)$ was determined in \cite[Theorem~6.8]{SCH}.}

For $r \in \mathbb{N}_0$ and $\alpha >0,$ define
\begin{align*}
W_c^{r+\alpha, \ast}(\mathbb{R}^+)&:=\,\left\{f \in X_c\::\: \sup_{v \in \mathbb{R}}|v^r[\Theta^\alpha_cf]^\wedge_{M_c}(c+iv)| < +\infty\right\} \\
&\,= \left\{f \in X_c\::\: \sup_{v \in \mathbb{R}}|v^{r+\alpha}[f]^\wedge_{M_c}(c+iv)| < +\infty\right\}.
\end{align*}
We note that $W_c^{r+\alpha,1}(\mathbb{R}^+) \subset W_c^{r+\alpha,\ast}(\mathbb{R}^+).$ Now the following proposition,
stated in \cite[Theorem~4.3]{BBMS4}, is a simple consequence of 
\cite[Theorem 5]{BBMS2}.
\begin{prop}\label{polyrate} (Polynomial rate).
Let $f \in \mathcal{M}^1_c.$ Then
$$f \in W_c^{r+\alpha, \ast}(\mathbb{R}^+) \:\Longleftrightarrow \: \dist_\infty(\Theta^\alpha_cf, B^1_{c,\sigma}) = \mathcal{O}(\sigma^{-r })
\qquad (\sigma \rightarrow +\infty).$$
\end{prop}

\begin{Remark}\label{rem_sob}
{\rm By inspection of its proof, it is easily seen that Theorem~\ref{estdistance} remains true if in its hypotheses
$W_c^{\alpha,1}(\mathbb{R}^+)$ is replaced by the larger space
$W_c^{\alpha,\ast}(\mathbb{R}^+)$. As a consequence, Lemma~\ref{lemma3} as well as Theorems~\ref{equiv2} and \ref{equiv1} hold with the
assumptions of Theorem~\ref{estdistance} being modified accordingly.}
\end{Remark}

\subsection{Characterization of the speed of convergence by function spaces}
\label{Sec:4.3}

Combining the results of Subsection~\ref{Sec:4.1} with those of Subsection~\ref{Sec:4.2}, we obtain a characterization of the speed of
convergence by function spaces.
\subsubsection{Case of the Mellin-even part}
The following statement is obtained by combining Proposition~\ref{evenodd}, Lemma~\ref{inverseformula} and Proposition~\ref{exactness}.
\begin{prop}\label{exactness3}
Let $\varphi \in \mathcal{K}_c(s) \cap \mathcal{M}^1_c.$ Then 
$$ \varphi_{c+}\in B^1_{0,2\pi\sigma} \: \Longleftrightarrow \: R_{c,\rho}[\varphi]=0 \; \hbox{ for all } \rho>\sigma.$$
\end{prop}
Next we characterize exponential rate of convergence.
\begin{prop}\label{hardy}
Let $\varphi \in \mathcal{K}_c(s) \cap \mathcal{M}^1_c.$ Then for $a >0,$ we have 
$$
\left.
\begin{array}{c}
f \in H^\ast_0(\mathbb{H}_a)\\[.8ex]
f(\cdot,0)\equiv \varphi_{c+}
\end{array} 
\right\}
\: \Longleftrightarrow \:R_{c,\sigma}[\varphi] = \mathcal{O}(e^{-2\pi a\sigma}) \quad (\sigma \rightarrow +\infty).$$
\end{prop}
{\bf Proof}. Let $f\in H^\ast_0(\mathbb{H}_a)$ such that $f(\cdot,0)\equiv \varphi_{c+}$. Then, by Proposition \ref{exprate}, 
we have $\dist_\infty(\varphi_{c+}, B^1_{0,2 \pi \sigma}) = \mathcal{O}(e^{-2\pi a\sigma})$ as $\sigma \rightarrow +\infty.$ 
Hence there exist $\sigma_0>0$ and $C>0$ such that 
$$\sup_{|v| \geq 2\pi\sigma}\left|[\varphi_{c+}]^\wedge_{M_0}(iv)\right| \leq C e^{-2\pi a\sigma} \qquad (\sigma\ge\sigma_0).$$
This implies that 
$$\left|[\varphi_{c+}]^\wedge_{M_0}(i2\pi k\sigma)\right| \leq C e^{-2\pi ak\sigma}$$
for all $k\in\mathbb{N}$ and $\sigma\ge\sigma_0$. Now, using (\ref{representation3}), we readily conclude that 
$R_{c,\sigma}[\varphi]=\mathcal{O}(e^{-2\pi a\sigma})$ as $\sigma\to +\infty.$

Conversely, if $R_{c,\sigma}[\varphi] = \mathcal{O}(e^{-2 \pi a\sigma})$ as $\sigma\to +\infty$, then there exist $\sigma_0>0$ and $C>0$
such that 
$$\left|R_{c,k\sigma}[\varphi]\right| \le C e^{-2\pi ak\sigma}\qquad (k\in\mathbb{N},\, \sigma\ge\sigma_0).$$
Therefore Lemma~\ref{inverseformula} applies, and (\ref{inverse2}) yields 
$$ [\varphi_{c+}]^\wedge_{M_0}(i2\pi\sigma)\,=\,\mathcal{O}(e^{-2\pi a\sigma}) \qquad (\sigma\to +\infty).$$
Hence $\dist_\infty(\varphi_{c+}, B^1_{0,\sigma})=\mathcal{O}(e^{-a\sigma})$ as $\sigma\to +\infty$.
Now Proposition~\ref{exprate} yields the desired assertion. $\:\Box$
\vskip0,4cm
Proposition~\ref{hardy} completes and generalizes a result in \cite[Theorem 6.5]{SCH}. The latter theorem needed the restriction 
$a \in ]0, 2\pi^2]$ in order to avoid the use of Riemann surfaces. With our concept of polar-analytic functions  
this restriction on $a$ does not arise.

Combining Theorem~\ref{equiv2} with Propositions~\ref{evenodd} and \ref{polyrate}, and taking Remark~\ref{rem_sob} into account,
we arrive at the following characterization of polynomial rate of convergence.

\begin{prop}\label{polyrate3} 
Let $\varphi \in \mathcal{K}_c(s) \cap \mathcal{M}^1_c.$ Then, for $r\in\mathbb{N}_0$ and $\alpha>1$, we have 
$$ \varphi_{c+}\in W_0^{r+\alpha,\ast}(\mathbb{R}^+) \: \Longleftrightarrow \: 
R_{c,\sigma}[\varphi]=\mathcal{O}(\sigma^{-r-\alpha})\quad (\sigma\to +\infty).$$
\end{prop}

\subsubsection{Case of the translated function $\tau^c_hf$ }
The propositions of Subsection~\ref{Sec:4.2} in conjunction with the results in \S\,\ref{subsub4.1.2} lead us to the following characterizations:
\begin{prop}\label{exactnees4}
Let $\varphi \in \mathcal{K}_c(s) \cap\mathcal{M}^1_c.$ Then we have
$$ \varphi\in B^1_{c,2\pi\sigma} \: \Longleftrightarrow \:
\left\{\begin{array}{l}
	R_{c,\rho}\left[\tau_h^c\varphi\right]=0\\[.8ex]
	\hbox{for all } \rho>\sigma, \, h\in [e^{-1/(2\rho)}, e^{1/(2\rho)}].
	\end{array} \right.
$$
\end{prop}
\smallskip

\begin{prop}\label{exprate4}
Let $\varphi \in \mathcal{K}_c(s) \cap\mathcal{M}^1_c.$ Then we have
$$
\left.\begin{array}{l}
	f \in H^\ast_c(\mathbb{H}_a)\\[.8ex]
	f(\cdot,0)\equiv \varphi
	\end{array}\right\}  \: \Longleftrightarrow \:
\left\{\begin{array}{l}
	R_{c,\sigma}\left[\tau_h^c\varphi\right]= \mathcal{O}(e^{-2\pi a\sigma}) \quad (\sigma\to +\infty)\\[.8ex]
	\hbox{uniformly for }  h\in [e^{-1/(2\sigma)}, e^{1/(2\sigma)}].
	\end{array} \right.
$$
\end{prop}
\smallskip

\begin{prop}\label{polyrate4}
Let $\varphi \in \mathcal{K}_c(s) \cap\mathcal{M}^1_c.$ Then, for $r\in\mathbb{N}_0$ and $\alpha>1$, we have
$$ \varphi\in W_c^{r+\alpha, \ast}(\mathbb{R}^+) \: \Longleftrightarrow \:
\left\{\begin{array}{l}
	R_{c,\sigma}\left[\tau_h^c\varphi\right]= \mathcal{O}(\sigma^{-r-\alpha}) \quad (\sigma\to +\infty)\\[.8ex]
	\hbox{uniformly for }  h\in [e^{-1/(2\sigma)}, e^{1/(2\sigma)}].
	\end{array} \right.
$$
\end{prop}
\section{Numerical examples}
\label{Sec:5}

In this section we give some examples illustrating the theory developed in previous sections.
Computations were performed with the help of Maple~16. 

\subsubsection*{Example 1: A Mellin-bandlimited integrand}
For $m\in\mathbb{N}$, we consider
$$ f_{2m}(r)\,:= 
\left\{
\begin{array}{cl} 
\displaystyle\left(\frac{\sin(\pi\log r)}{\pi\log r}\right)^{2m} &
\hbox{ \, if } r\ne 1,\\[2ex]
1 & \hbox{ \, if } r=1.
\end{array}\right.
$$
We want to compute the integral
$$ I_{2m}\,:=\, \int_0^\infty f_{2m}(r)\frac{dr}{r}\,.$$
Its exact value is given by (see \cite[p.~494, \S\,3.836/2]{GR})
$$ I_{2m}\,=\, \frac{m}{2^{2m-3}} \sum_{j=0}^{m-1} (-1)^j
\frac{(2m-2j)^{2m-2}}{j!(2m-1-j)!}\,.$$
In particular,
$$ I_2=1, \quad I_4=\frac{2}{3}, \quad I_6=\frac{11}{20}, \quad I_8=
\frac{151}{315}\,.$$
Note that $f_{2m}$ is a $0$-Mellin-even function which is Mellin-bandlimited
to $[-2\pi m, 2\pi m].$ The quadrature formula (\ref{ampsf}) now yields
\begin{equation}\label{qf}
I_{2m}\,=\, \frac{1}{\sigma} \left(1+2\sum_{k=1}^\infty f_{2m}\left(e^{k/\sigma}\right)\right)+
R_{0,\sigma}[f_{2m}]
\end{equation}
for any positive $\sigma$. According to the discussion in 
Subsection~\ref{sec:2.1}, we have
$R_{0,\sigma}[f_{2m}]=0$ as soon as $\sigma\ge m$. However, in computations we 
have to truncate the series in (\ref{qf}), and so there is always a
truncation error
$$T_{m,\sigma,K}\,:=\, 
\frac{2}{\sigma} \sum_{k=K}^\infty f_{2m}\left(e^{k/\sigma}\right).$$
By standard estimates and the integral comparison method, we find that
$$T_{m,\sigma,K}\,\le\, \frac{2}{\sigma}\sum_{k=K+1}^\infty
\left(\frac{\sigma}{\pi k}\right)^{2m}
\le\frac{2\sigma^{2m-1}}{\pi^{2m}} \int_K^\infty x^{-2m}dx=
\frac{2}{(2m-1)\pi}\left(\frac{\sigma}{K\pi}\right)^{2m-1}.
$$
Thus, in order to guarantee that the truncation error does not exceed
$10^{-\ell}$, say, we  should choose
$$ K \ge \frac{\sigma}{\pi} \left(\frac{2\cdot 10^\ell}{(2m-1)\pi}\right)^{1/(2m-1)}.$$
The total error is given by
$$ E_{m.\sigma,K}\,:=\, I_{2m} - \frac{1}{\sigma}\left(1 +2 \sum_{k=1}^K
f_{2m}\left(e^{k/\sigma}\right)\right).$$
For $m=4$ and $\ell=12$ computations provided the results shown
in Table~\ref{ex1}.
It is clearly seen that the total error reduces to the truncation error as
soon as $\sigma\ge 4.$

\begin{table}
{\footnotesize
\caption{Example~1 for $m=4$}\label{ex1}
\begin{tabular}{crr}\hline
\multicolumn{1}{c}{$\sigma$} & \multicolumn{1}{c}{$K$}&
\multicolumn{1}{c}{$E_{4,\sigma,K}$}\\\hline
0.50 & 6 & $-1.520635e+00$ \\
1.00& 12& $-5.206349e-01$ \\
1.50& 18& $-1.884456e-01$ \\
2.00& 24& $-4.761905e-02$ \\
2.50& 30& $-6.755332e-03$ \\
3.00& 36& $-3.968254e-04$ \\
3.50& 41& $-3.100198e-06$ \\
4.00& 47& $2.329011e-13$ \\
4.50& 53& $2.267646e-13$ \\
5.00& 59& $2.266369e-13$ \\
5.50& 65& $2.266218e-13$ \\
6.00& 71& $2.266187e-13$ \\
6.50& 77& $2.266178e-13$ \\
7.00& 82& $2.267143e-13$ \\
7.50& 88& $2.266709e-13$ \\
8.00& 94& $2.266480e-13$ \\\hline
\end{tabular}
}
\end{table}

\subsubsection*{Example 2: Integrand with a branch point}
We modify the function $f_8$ of the previous example by multiplying it
with a square root, which creates a branch point. More precisely, for
$a>0$, we introduce
$$ g_a(r)\,:=\, f_8(r)\,\sqrt{a^2 + (\log r)^2}.$$
We want to compute the integral
$$ J_a\,:=\, \int_0^\infty g_a(r)\,\frac{dr}{r}\,.$$
First we look for a polar-analytic extension of $g_a$. It is easy to see
that
$$ \widetilde{f}_8(r, \theta)\,:=\left\{
\begin{array}{cl} 
\displaystyle\left(\frac{\sin(\pi(\log r +i\theta))}{\pi(\log r+ i\theta)}\right)^8 
& \hbox{ \, if } (r,\theta)\in\mathbb{H}\setminus\{(1,0)\},\\[2ex]
1 & \hbox{ \, if } (r,\theta)=(1,0)
\end{array}\right.
$$
is a polar-analytic extension of $f_8$ to $\mathbb{H}.$
Therefore $\widetilde{g}_a$, defined by
$$ \widetilde{g}_a(r, \theta)\,:=\, \widetilde{f}_8(r,\theta)\, \sqrt{a^2
+(\log r + i\theta)^2},$$
is the right candidate for a polar-analytic extension of $g_a$. The
expression under the root vanishes if and only if $(r,\theta)=(1,\pm a)$.
Hence $\mathbb{H}_a$ is the largest strip in $\mathbb{H}$ on which
$\widetilde{g}_a$ is polar-analytic. Now it is easily verified that
$\widetilde{g}_a$ belongs to the Mellin--Hardy space $H^1_0(\mathbb{H}_a)$,
which is a subspace of $H^\ast_0(\mathbb{H}_a)$, but it does not belong to
$H^\ast_0(\mathbb{H}_b)$ for $b>a$.

Since $g_a\in X_0$, we have $M^\ast_0[g_a](it)\equiv M_0[g_a](it)$. This
allows us to conclude with the help of Theorem~\ref{expopw} that
$g_a\in \mathcal{K}_0(s)\cap \mathcal{M}_0^1$ for any positive $s$.
Hence the theory established in Sections~\ref{Sec:3}--\ref{Sec:4} is applicable to
$g_a$.

We note that $g_a$ is again a $0$-Mellin-even function. By the quadrature
formula (\ref{ampsf}) we have
$$ J_a\,=\, \frac{1}{\sigma}\left[a + 2 \sum_{k=1}^\infty
f_8(e^{k/\sigma}) \sqrt{a^2 + (k/\sigma)^2}\right] + R_{0,\sigma}[g_a].$$
The exact value of $J_a$ is not known. For numerically given $a$, we
used Maple for gaining $J_a$ up to $40$ decimal places. The total error in our
computation by formula (\ref{ampsf}) will be
$$ E_{a,\sigma,K} \,:=\, J_a - \frac{1}{\sigma}\left[a + 2 \sum_{k=1}^K
f_8(e^{k/\sigma}) \sqrt{a^2+(k/\sigma)^2}\right]
\,=\, R_{0,\sigma}[g_a] + T_{a,\sigma,K},$$
where
$$ T_{a,\sigma,K}\,:=\, \frac{2}{\sigma} \sum_{k=K+1}^\infty
f_8(e^{k/\sigma}) \sqrt{a^2+ (k/\sigma)^2}$$
is the truncation error. Since
$$ f_8(e^{k/\sigma})\,\le\, \frac{1}{(\pi k/\sigma)^8} \quad \hbox{ and }
\quad \sqrt{a^2+ (k/\sigma)^2}\,\le\, \frac{k}{\sigma} +
\frac{a^2\sigma}{2k}\,,$$
we find by using the integral comparison method that
\begin{align*}
T_{a,\sigma,K} & \le\, \frac{2}{\pi^8}
\int_K^\infty\left[\left(\frac{\sigma}{x}\right)^7 +
\frac{a^2}{2}\left(\frac{\sigma}{x}\right)^9\right]
\frac{dx}{\sigma}\\[2ex]
&=\,\frac{1}{\pi^8}\left[\frac{1}{3}\left(\frac{\sigma}{K}\right)^6 + 
\frac{a^2}{8} \left(\frac{\sigma}{K}\right)^8\right].
\end{align*}

Proposition~\ref{hardy} tells us that $R_{0,\sigma}[g_a]=
\mathcal{O}(e^{-2\pi a\sigma})$ as $\sigma\to +\infty$. Hence, in order that
the truncation error does not exceed the remainder asymptotically, we
choose $K$ such that
$$
\frac{1}{3}\left(\frac{\sigma}{K}\right)^6 + \frac{a^2}{8}
\left(\frac{\sigma}{K}\right)^8 \,\le\, \frac{\pi^8}{10} \, e^{-2\pi
a\sigma}.$$

The asymptotic result for the remainder and the choice of $K$ for
controlling the truncation error suggest that
$$ E_{a,\sigma,K}\, \approx\, C\,e^{-2\pi a\sigma}$$
for large $\sigma$. In our numerical experiments we check this behavior in
two ways by computing
\begin{equation}\label{test}
C\,:=\,E_{a,\sigma,K}\,e^{2\pi a\sigma} \quad \hbox{ and }\quad
\hbox{ rate} \,:=\, -\frac{\log |E_{a,\sigma,K}|}{\sigma}\,.
\end{equation}
While $C$ should remain bounded, the rate should approach $2\pi a$ as
$\sigma\to +\infty$. As an immediate consequence of (\ref{test}), we have
$$ \hbox{rate } =\, 2\pi a - \frac{\log |C|}{\sigma}\,.$$
Hence for bounded $\sigma$, say $2\le\sigma\le 15$, the numbers rate can
be close to $2\pi a$ only if $|C|$ is close to $1$. This explains the
behavior of the numbers in the last two columns of
Tables~\ref{ex2.1}--\ref{ex2.3}. In view of Theorem~\ref{expopw}, we
expect that $C$ depends on $\|g_a\|_{H_0^\ast(\mathbb{H}_a)}.$ It seems
that this expression is growing considerably when $a$ moves from $1/2$
to $1$.

\begin{table}[h]
{\footnotesize
\caption{Example~2 for $a=\frac{1}{2}$, rate $\to \pi=3.141592\dots$}\label{ex2.1}
\begin{tabular}{rrrrr}
\hline
&&&&\\[-1.8ex]
\multicolumn{5}{c}{%
$J_a= 0.2552373684721620868389158816136888733878$}
\\[.2ex]\hline
\multicolumn{1}{r}{$\sigma$} & \multicolumn{1}{r}{$K$} &
\multicolumn{1}{c}{$E_{a,\sigma,K}$} & \multicolumn{1}{c}{$C$} &
\multicolumn{1}{c}{rate}\\\hline
$2$ & $8$ & $-1.385e-02$ & $-7.415e+00$ & $2.139842$ \\
$3$ & $4$ & $4.441e-04$ & $5.503e+00$ & $2.573168$ \\
$4$ & $9$ & $6.830e-06$ & $1.959e+00$ & $2.973533$ \\
$5$ & $19$ & $1.453e-07$ & $9.641e-01$ & $3.148914$ \\
$6$ & $37$ & $4.117e-09$ & $6.321e-01$ & $3.218038$ \\
$7$ & $73$ & $1.274e-10$ & $4.526e-01$ & $3.254852$ \\
$8$ & $141$ & $4.216e-12$ & $3.467e-01$ & $3.274013$ \\
$9$ & $267$ & $1.478e-13$ & $2.812e-01$ & $3.282555$ \\
$10$ & $500$ & $5.372e-15$ & $2.365e-01$ & $3.285766$ \\
$11$ & $927$ & $1.995e-16$ & $2.033e-01$ & $3.286429$ \\
$12$ & $1707$ & $7.525e-18$ & $1.774e-01$ & $3.285687$ \\
$13$ & $3122$ & $2.883e-19$ & $1.573e-01$ & $3.283862$ \\
$14$ & $5675$ & $1.120e-20$ & $1.414e-01$ & $3.281298$ \\
$15$ & $10264$ & $4.400e-22$ & $1.285e-01$ & $3.278358$ \\
\hline
\end{tabular}
}
\end{table}

\begin{table}[h]
{\footnotesize
\caption{Example~2 for $a=\frac{5}{8}$, rate $\to \frac{5\pi}{4}=3.926990\dots$}\label{ex2.2}
\begin{tabular}{rrrrr}
\hline
&&&&\\[-1.8ex]
\multicolumn{5}{c}{$J_a= 0.3123770437749010235851625171708586776416$}
\\[.2ex]\hline
\multicolumn{1}{r}{$\sigma$} & \multicolumn{1}{r}{$K$} &
\multicolumn{1}{c}{$E_{a,\sigma,K}$} & \multicolumn{1}{c}{$C$} &
\multicolumn{1}{c}{rate}\\\hline
$2$ & $3$ & $-2.172e-02$ & $-5.596e+01$ & $1.914666$ \\
$3$ & $6$ & $2.043e-04$ & $2.671e+01$ & $2.831941$ \\
$4$ & $15$ & $1.962e-06$ & $1.302e+01$ & $3.285439$ \\
$5$ & $36$ & $1.730e-08$ & $5.827e+00$ & $3.574493$ \\
$6$ & $81$ & $2.082e-10$ & $3.558e+00$ & $3.715444$ \\
$7$ & $182$ & $2.856e-12$ & $2.478e+00$ & $3.797346$ \\
$8$ & $400$ & $4.219e-14$ & $1.858e+00$ & $3.849564$ \\
$9$ & $865$ & $6.546e-16$ & $1.463e+00$ & $3.884719$ \\
$10$ & $1849$ & $1.051e-17$ & $1.192e+00$ & $3.909403$ \\
$11$ & $3912$ & $1.733e-19$ & $9.976e-01$ & $3.927212$ \\
$12$ & $8212$ & $2.916e-21$ & $8.520e-01$ & $3.940341$ \\
\hline
\end{tabular}
}
\end{table}

\begin{table}[h]
{\footnotesize
\caption{Example~2 for $a=1$, rate $\to 2\pi=6.283185\dots$}\label{ex2.3}
\begin{tabular}{rrrrr}
\hline
&&&&\\[-1.8ex]
\multicolumn{5}{c}{$J_a = 0.4876105654991947134580915823151850342698$}
\\[.2ex]\hline
\multicolumn{1}{r}{$\sigma$} & \multicolumn{1}{r}{$K$} &
\multicolumn{1}{c}{$E_{a,\sigma,K}$} & \multicolumn{1}{c}{$C$} &
\multicolumn{1}{c}{rate}\\\hline
$2$ & $5$ & $-4.256e-02$ & $-1.220e+04$ & $1.578402$ \\
$3$ & $19$ & $-1.580e-04$ & $-2.426e+04$ & $2.917691$ \\
$4$ & $71$ & $1.248e-07$ & $1.026e+04$ & $3.974063$ \\
$5$ & $250$ & $7.809e-11$ & $3.438e+03$ & $4.654642$ \\
$6$ & $854$ & $8.050e-14$ & $1.898e+03$ & $5.025082$ \\
$7$ & $2838$ & $9.886e-17$ & $1.248e+03$ & $5.264689$ \\
$8$ & $9241$ & $1.333e-19$ & $9.014e+02$ & $5.432688$ \\
$9$ & $29625$ & $1.907e-22$ & $6.905e+02$ & $5.556805$ \\
$10$ & $93800$ & $2.841e-25$ & $5.509e+02$ & $5.652030$ \\
\hline
\end{tabular}
}
\end{table}

\subsubsection*{Example 3: An integral representing the gamma function}
For  complex arguments $z$ with $\hbox{Re}\,z >0$ the gamma function 
was defined by Euler (see, e.g., \cite[p.~942, \S\,8.310]{GR}) as
\begin{equation}\label{gamma}
\Gamma (z) = \int_0^\infty e^{-r}r^{z-1}dr.
\end{equation}
Writing $f(r):=e^{-r}$ and $z=c+it$ with $c>0, \,t\in\mathbb{R}$, we see
that
$\Gamma (z) = [f]_{M_c}^\wedge(c+it)$. 
We note that $\widetilde{f}$, defined by
$$ \widetilde{f}(r,\theta)\,:=\, \exp\left(-r e^{i\theta}\right)\,=\,
e^{-r\cos \theta} e^{-ir \sin \theta}, \qquad (r,\theta)\in\mathbb{H},$$
is a polar-analytic extension of $f$. It can be easily verified that for
any $a\in ]0, \pi/2[$, we have $\widetilde{f}\in H_c^1(\mathbb{H}_a)$ but
$\widetilde{f}\not\in H_c^\ast(\mathbb{H}_{\pi/2})$ since condition (d) of
Definition~\ref{intermediate} fails. 

As in the previous example, we can conclude with the help of
Theorem~\ref{expopw} that $f\in \mathcal{K}_c(s)\cap \mathcal{M}_c^1$ for
every positive $s$. Now Proposition~\ref{exprate4} yields that
$$R_{c,\sigma}[f] = \mathcal{O}(e^{-2\pi a \sigma})\qquad  (\sigma \rightarrow +\infty)$$
for each $a \in ]0, \pi/2[.$
For $c = 1/2$ the quadrature formula (\ref{ampsf}) reads as
\begin{equation}\label{quad_gamma}
\int_0^\infty f(r) r^{1/2}\, \frac{dr}{r}\, = \, 
\frac{1}{\sigma}\sum_{k=-\infty}^\infty \exp\left(-e^{k/\sigma}\right) e^{k/(2\sigma)} + 
R_{1/2,\sigma}[f].
\end{equation}
By (\ref{gamma}), the left-hand side is equal to
$\Gamma(1/2)=\sqrt{\pi}$. This time the integrand is not Mellin-even.
Therefore we prefer an appropriate asymmetric trunction of the series in
(\ref{quad_gamma}), which gives a total error
$$E_{\sigma, N,K} := \sqrt{\pi} - \frac{1}{\sigma}\sum_{k=-N}^K 
\exp\left(-e^{k/\sigma}\right) e^{k/(2 \sigma)}.$$

By standard estimates using the integral comparison method, we find that
the choice
$$ N := \left\lceil 2\pi^2\sigma^2 +2\sigma \log \frac{15}{4}\right\rceil,
\quad K:= \left\lceil \sigma\,\log\frac{N}{2\sigma}\right\rceil$$
will produce a truncation error that does not exceed the remainder 
asymptotically. Here $\lceil\cdot\rceil$ denotes the ceiling function
mapping $x$ to the least integer that is greater than or equal to $x$.
Analogously to the previous example, we also compute
$$ C\,:=\, E_{\sigma,N,K}\,e^{\pi^2 \sigma} \quad \hbox{ and }\quad
\hbox{rate}\, :=-\,\frac{\log |E_{\sigma,N,K}|}{\sigma}\,.$$
By our theory, the numbers rate should become larger than $2\pi a$ for
any $a< \pi/2$ as $\sigma\to +\infty.$ 
On the other hand, it can be shown that the $1/2$-Mellin-even part
$f_{1/2+}$ does not have an extension belonging to
$H_0^\ast(\mathbb{H}_{\pi/2})$. Therefore Propositions~\ref{hardy} and
\ref{exprate4} imply that the numbers rate must converge to $\pi^2$ as
$\sigma\to +\infty.$

Table~\ref{ex3.1} shows that for small values of $\sigma$ and,
consequently, with relatively short sums, we obtain already satisfactory
approximations to the integral. However, the numbers rate are still
considerably smaller than $\pi^2$. Therefore we continued with
larger values of $\sigma$, setting Maple's environment variable
{\it Digits}\,$:=80$ 
in order to suppress round-off errors. Table~\ref{ex3.2} shows
that the approximations to the integral become excellent and
the numbers rate get much closer to $\pi^2$.

\begin{table}[h]
{\footnotesize
\caption{Example~3 for small $\sigma$, rate $< \pi^2=9.869604\dots$}\label{ex3.1}
\begin{tabular}{rrrrrr}
\hline
\multicolumn{1}{c}{$\sigma$} & \multicolumn{1}{r}{$N$} &
\multicolumn{1}{r}{$K$} &
\multicolumn{1}{c}{$E_{\sigma,N,K}$} & \multicolumn{1}{c}{$C$} &
\multicolumn{1}{c}{rate}\\\hline
$0.25$ & $2$ & $1$ & $-3.038e-01$ & $-3.583$ & $4.765302$ \\
$0.50$ & $7$ & $1$ & $-3.125e-02$ & $-4.346$ & $6.931166$ \\
$0.75$ & $14$ & $2$ & $2.746e-03$ & $4.502$ & $7.863614$ \\
$1.00$ & $23$ & $3$ & $-1.219e-04$ & $-2.356$ & $9.012572$ \\
$1.25$ & $35$ & $4$ & $1.160e-05$ & $2.644$ & $9.091811$ \\
$1.50$ & $49$ & $5$ & $-1.103e-06$ & $-2.966$ & $9.144985$ \\
$1.75$ & $66$ & $6$ & $1.600e-07$ & $5.090$ & $8.941767$ \\
$2.00$ & $85$ & $7$ & $-1.100e-08$ & $-4.200$ & $9.162685$  \\
\hline
\end{tabular}
}
\end{table}

\begin{table}[h]
{\footnotesize
\caption{Example~3 continued, rate $\to \pi^2=9.869604\dots$}\label{ex3.2}
\begin{tabular}{rrrrrr}
\hline
\multicolumn{1}{r}{$\sigma$} & \multicolumn{1}{r}{$N$} &
\multicolumn{1}{r}{$K$} &
\multicolumn{1}{c}{$E_{\sigma,N,K}$} & \multicolumn{1}{c}{$C$} &
\multicolumn{1}{c}{rate}\\\hline
$2$ & $85$ & $7$ & $-1.135e-08$ & $-4.242$ & $9.147070 $ \\
$3$ & $186$ & $11$ & $-2.546e-13$ & $-1.840$ & $9.666308 $ \\
$4$ & $327$ & $15$ & $-2.862e-17$ & $-3.999$ & $9.523119 $ \\
$5$ & $507$ & $20$ & $-1.539e-22$ & $-0.416$ & $10.045137 $ \\
$6$ & $727$ & $25$ & $-1.680e-26$ & $-0.877$ & $9.891425 $ \\
$7$ & $986$ & $30$ & $4.948e-30$ & $4.996$ & $9.639809 $ \\
$8$ & $1285$ & $36$ & $1.358e-34$ & $2.651$ & $9.747734 $ \\
$9$ & $1623$ & $41$ & $5.354e-39$ & $2.021$ & $9.791439 $ \\
$10$ & $2001$ & $47$ & $5.676e-43$ & $4.142$ & $9.727498 $ \\
$11$ & $2418$ & $52$ & $2.957e-47$ & $4.171$ & $9.739769 $ \\
$12$ & $2875$ & $58$ & $-1.102e-51$ & $-3.007$ & $9.777868 $ \\
$13$ & $3371$ & $64$ & $-7.861e-57$ & $-0.415$ & $9.937342 $ \\
$14$ & $3906$ & $70$ & $2.173e-60$ & $2.215$ & $9.812793 $ \\
$15$ & $4481$ & $76$ & $-9.514e-65$ & $-1.875$ & $9.827685 $ \\
\hline
\end{tabular}
}
\end{table}

\subsubsection*{Example 4: Integrand belonging to a Mellin--Sobolev space}
Consider the function $g:\mathbb{R}^+ \rightarrow \mathbb{R}$ 
defined by
\begin{eqnarray*}
g(r) = \left\{\begin{array}{ll} r\log^2r, &~0<r<1,
 \\[1ex] r^{-1} \log^2r, &~r\geq 1.
\end{array} \right.
\end{eqnarray*}
It is seen to be $0$-Mellin-even. By a straightforward calculation, we find
that
$$[g]^\wedge_{M_0}(iv) \,=\, 4\,\frac{1-3v^2}{(1+v^2)^3} \qquad 
(v \in \mathbb{R}),$$
and so
$$ \int_0^\infty g(r)\,\frac{dr}{r}\,=\,[g]_{M_0}^\wedge(0)\,=\,4.$$
Since we know $g$ and its Mellin transform explicitly, we can readily
verify that $g\in \mathcal{K}_0(s)\cap \mathcal{M}_0^1$ for any positive
$s$. Furthermore, it can be shown that
$g\in W_0^{4,1}(\mathbb{R}^+) \subset W_0^{4,\ast}(\mathbb{R}^+)$
but $g\not\in W_0^{4+\alpha, \ast}(\mathbb{R}^+)$ for  $\alpha>0$.

The remainder of the quadrature formula (\ref{ampsf}) can now be written as
\begin{align*}
 R_{0,\sigma}[g] &=\,\int_0^\infty g(r)\,\frac{dr}{r} -\frac{1}{\sigma}\sum_{k=-\infty}^\infty g(e^{k/\sigma})\\[1ex]
	&=\, 4 - \frac{2}{\sigma}\sum_{k=1}^\infty\bigg(\frac{k}{\sigma}\bigg)^2 e^{-k/\sigma} \\[1ex]
	&=\, 4 - \frac{2}{\sigma^3} \cdot \frac{e^{-2/\sigma} + e^{-1/\sigma}}{(1- e^{-1/\sigma})^3}\,.
\end{align*}
An expansion of the right-hand side yields
\begin{equation}\label{expansion}
R_{0,\sigma}[g]\,=\,\frac{1}{60\,\sigma^4} - \frac{1}{756\,\sigma^6} + \mathcal{O}\left(\sigma^{-8}\right)\qquad (\sigma\to +\infty).
\end{equation}
Next, from Theorem~\ref{estdistance}, we deduce that
\begin{align*}
\left|R_{0,\sigma}[g]\right| &\le\, \frac{2 \zeta(4)}{(2\pi\sigma)^4}\,\dist_\infty(\Theta_0^4 g, B^1_{0, 2\pi\sigma})\\[1ex]
	&=\,\frac{2 \zeta(4)}{(2\pi\sigma)^4}\, \sup_{|v|\ge 2\pi\sigma}\left|4v^4\,\frac{1-3v^2}{(1+v^2)^3}\right|\,=\,
\frac{1}{60\, \sigma^4}\,.
\end{align*}
Now comparison with (\ref{expansion}) shows that the estimate (\ref{Sob_bound}) is best possible in the sense described in
Theorem~\ref{estdistance}.

The precision of this estimate is  illustrated in Table~\ref{ex4a}. The third column shows that the upper bound in (\ref{Sob_bound}) becomes
very close to the true value of the remainder as $\sigma$ grows. Consequently, the factor of overestimation, defined by
$$ \hbox{overestimation}\,:=\,\frac{\hbox{upper bound}}{\hbox{remainder}}\,=\, \frac{\frac{1}{60\sigma^4}}{R_{0,\sigma}[g]}\,,$$
becomes very close to $1$. The fifth column shows $C:= R_{0,\sigma}[g]\,\sigma^4$, which converges rapidly to $1/60$ as $\sigma\to +\infty$.
%
%

\begin{table}[h]
{\footnotesize
\caption{Example~4: $C \to \frac{1}{60}=0.01666\dots$}\label{ex4a}
\begin{tabular}{rrrrr}
\hline
\multicolumn{1}{r}{$\sigma$} & \multicolumn{1}{c}{$R_{0,\sigma}[g]$} &
\multicolumn{1}{c}{upper bound (\ref{Sob_bound})} & \multicolumn{1}{c}{overestimation}
 &\multicolumn{1}{c}{$C$}
\\\hline
$2$ & $1.021267e-03$ & $1.041667e-03$ & $1.019974824393$ & $0.016340272591$ \\
$4$ & $6.478229e-05$ & $6.510417e-05$ & $1.004968650783$ & $0.016584265244$ \\
$8$ & $4.063969e-06$ & $4.069010e-06$ & $1.001240599981$ & $0.016646015620$ \\
$16$ & $2.542343e-07$ & $2.543132e-07$ & $1.000310052376$ & $0.016661500729$ \\
$32$ & $1.589334e-08$ & $1.589457e-08$ & $1.000077506994$ & $0.016665374984$ \\
$64$ & $9.933915e-10$ & $9.934107e-10$ & $1.000019376367$ & $0.016666343733$ \\
$128$ & $6.208787e-11$ & $6.208817e-11$ & $1.000004844068$ & $0.016666585933$ \\
$256$ & $3.880506e-12$ & $3.880511e-12$ & $1.000001211016$ & $0.016666646483$ \\
$512$ & $2.425318e-13$ & $2.425319e-13$ & $1.000000302754$ & $0.016666661621$ \\
$1024$ & $1.515824e-14$ & $1.515825e-14$ & $1.000000075688$ & $0.016666665405$ \\
$2048$ & $9.473903e-16$ & $9.473903e-16$ & $1.000000018922$ & $0.016666666351$ \\
$4096$ & $5.921189e-17$ & $5.921189e-17$ & $1.000000004731$ & $0.016666666588$ \\
$8192$ & $3.700743e-18$ & $3.700743e-18$ & $1.000000001183$ & $0.016666666647$ 
\\\hline
\end{tabular}
}
\end{table}

\section{A short biography of Helmut Brass 1936--2011}\label{sec_Brass}

Helmut Brass (originally written as Bra{\ss}) was born in Hannover, Germany, in 1936. After completing Oberrealschule, he worked in a firm that produced 
and recycled copper cables. Since already as a schoolboy he was very interested in chemistry, he then enrolled for this subject at the University of 
Hannover. But soon after he had started, he realized that his true talent was mathematics and turned to it. In 1962 he graduated with a diploma in 
mathematics and continued as a Scientific Assistant under the supervision of Wilhelm Quade. 
In 1965 Brass received his doctoral degree in mathematics with a thesis on approximation by a linear combination of projection operators. In 1968 he 
acquired Habilitation and became a University Dozent.

In 1970, Brass was appointed as a professor at the Technical University of Clausthal and in 1974 he followed the offer of a chair at the University of 
Osnabrück, but already in 1977 he accepted a chair at the University of Braunschweig where he stayed until his retirement in 2002.

The research field of Helmut Brass comprised interpolation and approximation with special emphasis on quadrature. In particular, he studied optimal and 
nearly optimal quadrature formulae for various classes of functions, properties of the remainder functional such as positivity and monotonicity, best 
or asymptotically best error estimates for classical quadrature formulae and exact rates of convergence under side conditions on the involved function 
such as periodicity, convexity or bounded variation. He published about 50 research papers and two distinguished books: Quadraturverfahren in 1977 
\cite{BRA2} and (with K.\ Petras) Quadrature Theory in 2011, a product of almost 20 years of joint work. 
Brass also edited two Proceedings of Oberwolfach Conferences on numerical integration. Furthermore he wrote a fascinating booklet with eleven lectures 
on Bernoulli polynomials designed for the training of students in \emph{Pro\-semi\-nars}.

Brass had twelve research students who graduated with a doctoral degree under his supervision. Four of them acquired the Habilitation degree and one 
was also awarded with the title of a University Professor.

In 1963, Brass married Gisela Lueder. They had studied together in Hannover. She was a Gymnasium teacher of mathematics and chemistry. They had two sons, 
Stefan and Peter, both now being professors of computer science, one in Halle (Germany), and the other in New York. 

In 2008 a stroke of fate met the whole family, when Mrs.~Gisela Brass died all of a sudden, a shock from which Helmut never recovered. 
He passed away in Halle on October~30, 2011 in the house of his elder son Stefan.
\medskip

I (G.\,S.) met Helmut Brass for the first time in Oberwolfach in 1977 at a conference on Numerical Methods in Approximation Theory and again in 
Oberwolfach in 1978, 1981, 1987, 1992 and 2001 at conferences on Numerical Integration, as well as on a few other occasions.

A few months preceding our first meeting in 1977, Brass' book \cite{BRA2} had appeared. I was impressed by its systematic composition and the 
wealth of results. 
There I found an open problem for which I had an idea. When I contacted Brass, we maintained scientific correspondence over a period of more than two 
years that resulted in two joint papers
\cite{BRASCH1}, \cite{BRASCH2}.

Later I profited from work of Brass \cite{BRA1} in my collaboration with Q.\,I.~Rahman, when we characterized the speed of convergence of the 
trapezoidal formula and related quadrature methods in terms of function spaces; see, e.g., \cite{RS2}.

In 1979 Brass visited me in Erlangen and gave a talk in our Mathematical Colloquium. He reciprocated by inviting me to a colloquium talk in 
Braunschweig in 1981.
\medskip

I (P.L.\,B.) invited Wilhelm Quade (1898--1975), a discoverer of splines, to my first conference at Oberwolfach (1963). He brought along with him his 
research assistant Helmut Brass, at the time a shy young man.
I also invited Brass on the recommendation of my colleague Rolf Nessel to my Oberwolfach conference in 1980.

The final time I met him was ca.~1988 when he invited me to give a colloquium lecture at the Technische Hochschule Braunschweig. 
I accepted his kind invitation although I had turned down all similar foregoing invitations in Germany for some
twenty years.

Applications of the uniform boundedness principle of functional analysis, a chief area of research  of
Lehrstuhl A f\"ur Mathematik in the eighties \cite{DNvW}, were motivated by three basic papers of Helmut Brass, 
in which he gave best possible error estimates for quadrature rules \cite{BRA0}, \cite{BRA}, \cite{BRA1}, as well as his book \cite{BRA2}.

{\bf Acknowledgements}. 
Carlo Bardaro and Ilaria Mantellini have been partially supported by the ``Gruppo Nazionale per l'Analisi Matematica e Applicazioni'' (GNAMPA) of the 
``Istituto Nazionale di Alta Matematica'' (INDAM) as well as by the Department of Mathematics and Computer Sciences of the University of Perugia.

\bibliographystyle{spmpsci}      

\begin{thebibliography}{99}
\bibitem{BBM2} C. Bardaro, P.L. Butzer and I. Mantellini, The foundation of the fractional calculus in Mellin transform setting with applications,  
J. Fourier Anal. Appl, 21, 961--1017 (2015)
\bibitem{BBMS} C. Bardaro, P.L. Butzer, I. Mantellini and G. Schmeisser, On the Paley-Wiener theorem in the Mellin transform setting, 
J. Approx. Theory, 207, 60--75 (2016)
\bibitem{BBMS2} C. Bardaro, P.L. Butzer, I. Mantellini and G. Schmeisser, Mellin analysis and its basic associated metric. Applications to sampling theory, 
Analysis Math., 42(4), 297--321 (2016)
\bibitem{BBMS3} C. Bardaro, P.L. Butzer, I. Mantellini and G. Schmeisser, A fresh approach to the Paley-Wiener theorem for Mellin transforms and the 
	Mellin-Hardy spaces, to appear in Math. Nachr. (2017), DOI: 10.1002/mana.201700043
\bibitem{BBMS4} C. Bardaro, P.L. Butzer, I. Mantellini and G. Schmeisser, A generalization of the Paley--Wiener theorem for Mellin transforms and 
	metric characterization of function spaces, Frac. Calc. Appl. Anal., 20(5), 1216--1238 (2017)
\bibitem{BDR} C. de Boor, R.A. DeVore and A. Ron, Approximation from shift-invariant subspaces of $L_2(\mathbb{R}^d)$, Trans. Amer. Math. Soc. 341(2),
	787--806 (1994)
\bibitem{BRA0} H. Brass, Eine Fehlerabsch\"{a}tzung zum Quadraturverfahren von Clenshaw   und Curtis, Numer. Math. 21(5),  397--403 (1973)
\bibitem{BRA2} H. Brass, Quadraturverfahren, Vandenhoeck \& Ruprecht, G\"{o}ttingen, (1977)
\bibitem{BRA} H. Brass, Der Wertebereich des Trapezverfahrens, in ``G. H\"{a}mmerlin (Ed.), Numerische Integration,  Proc. Conf. Math. Res. Inst., Oberwolfach, Germany, Oct. 1-7, 1978, Birkh\"{a}user Verlag, Basel-Boston, Mass., 98--108 (1979)
\bibitem{BRA1} H. Brass, Umkehrs\"{a}tze beim Trapezverfahren, Aequationes Math. 18, 338--344 (1978)
\bibitem{BRASCH1} H. Brass and G. Schmeisser, The definiteness of Filippi's quadrature formulae  and related problems, in G.~H\"{a}mmerlin (Ed.), Numerische Integration,
  Proc. Conf. Math. Res. Inst., Oberwolfach, Germany, Oct. 1--7,
  1978, Vol. 45 of Internat. Ser. Numer. Math., Birkh\"{a}user, Basel-Boston,
  Mass.,  109--119 (1979)
 \bibitem{BRASCH2} H. Brass and G. Schmeisser, Error estimates for interpolatory quadrature formulae,   Numer. Math., 37(3),  371--386 (1981)
\bibitem{BDFHSS} P.L. Butzer, M. Dodson, P. Ferreira, J. Higgings, G. Schmeisser and R.L. Stens, Seven pivotal theorems of Fourier analysis, signal analysis, numerical analysis and number theory: their interconnections. Bull. Math. Sci.,  4(3), 481-525 (2014)
\bibitem{BJ1} P.L. Butzer and S. Jansche,  A direct approach to the Mellin transform, J. Fourier Anal. Appl., 3, 325--375 (1997)
\bibitem{BJ2} P.L. Butzer and S. Jansche, The finite Mellin transform, Mellin-Fourier series, and the Mellin-Poisson summation formula, Rend. Circ. Mat. Palermo, Serie II, Suppl. Vol. 52, 55--81 (1998)
\bibitem{BJ4} P.L. Butzer and S. Jansche, A self-contained approach to Mellin transform analysis for square integrable functions, applications, 
Integral Transforms Spec. Funct., 8, 175--198  (1999)
\bibitem{BN} P.L. Butzer and R.J. Nessel, Fourier analysis and approximation I, Birkh\"{a}user   Verlag, Basel; Academic Press, New York  (1971)
\bibitem{BSS-Marvasti} P.L. Butzer, G. Schmeisser and R.L. Stens, An introduction to sampling analysis, in F. Marvasti (ed.) ``Nonuniform Sampling: Theory and Practice'',
	(ISBN 0-306-46445-4), Kulwer Academic/Plenum Publishers, New York, 17--121 (2001)
\bibitem{BSS0} P.L. Butzer, G. Schmeisser and  R.L. Stens, Basic relations valid for the
  Bernstein space $B^p_{\sigma}$ and their extensions to functions from
  larger spaces with error estimates in terms of their distances from
  $B^p_{\sigma}$, J. Fourier Anal. Appl., 19(2), 333--375  (2013)
\bibitem{BSS} P.L. Butzer, G. Schmeisser and R.L. Stens, The distance between the general Poisson summation formula and that for bandlimited functions; 
applications to quadrature formulae, Appl. Comput. Harmonic. Anal. (2017), DOI: 10.1016/j.acha.2017.02.001
\bibitem{BW} P.L. Butzer and U. Westphal, An access to fractional differentiation via fractional differece quotients, in "Fractional calculus and its Applications", Proc. conf. New Haven, Lecture Notes in Math, 457, 116--145, Springer, Heidelberg  (1975)
\bibitem{DL} R.A. DeVore and  G.G. Lorentz, Constructive Approximation, Springer-Verlag,
 Berlin (1993)
\bibitem{DNvW}
W.~Dickmeis, R.~J. Nessel and E.~van Wickeren, Quantitative extensions of the
  uniform boundedness principle, Jahresber. Deutsch. Math.-Verein. 89(3), 105--134 (1987)
\bibitem{DRS} D.P. Dryanov, Q.I. Rahman and G. Schmeisser, Converse theorems in the theory of approximate integration, Constructive Approx., 6,  321--334 (1990)
\bibitem{Gautschi} W. Gautschi, Quadrature formulae on half-infinite intervals, BIT, 31, 438--446 (1991)
\bibitem{GPS} H-J. Glaeske, A.P. Prudnikov and K.A. Skornik, Operational calculus
and related topics, Chapman and Hall, CRC, Boca Raton, FL (2006)
\bibitem{GR} I.S. Gradshteyn and I.M. Ryzhik, Table of Integrals, Series and Products, Fifth Edition, Academic Press, New York, London 1994
\bibitem{LY} J.N. Lyness, The calculation of Fourier coefficients by the M\"{o}bius inversion of the Poisson summation formula. I. Functions whose early derivatives are continuous, Math. Comp. 24, 101--135  (1970)
\bibitem{LS} J.H. Loxton and J.W. Sanders, The kernel of a rule of approximate integration, J. Austral. Math. Soc. Ser. B, 21(3),  257--267 (1980)
\bibitem{LS1} J.H. Loxton and J.W. Sanders, On an inversion theorem of {M}\"obius, J. Austral.   Math. Soc. Ser. A, 30(1),  15--32  (1980)
\bibitem{MA} R.G. Mamedov, The Mellin Transform and Approximation Theory, (in Russian), "Elm", Baku, (1991)
\bibitem{MM} G. Mastroianni and G. Monegato, Truncated quadrature rules over $(0, \infty)$ and Nystr\"{o}m-type methods, SIAM J. Numer. Anal., 41(5), 1870--1892 (2003)
\bibitem{MNM} G. Mastroianni, I. Notarangelo and G.V. Milovanovi\'c, Gaussian quadrature rules with an exponential weight on the real semiaxis,
	IMA J. Numer. Anal., 34, 1654--1685 (2014)
\bibitem{RS} Q.I. Rahman and G. Schmeisser, On a Gaussian quadrature formula for entire functions of exponential type, in ``Numerical Methods of Approximation Theory'', Vol. 8, ISNM81, L. Collatz et al. (eds), Birkhauser-Verlag, Basel, 155--168 (1987)
\bibitem{RS2} Q.I. Rahman and G. Schmeisser, Characterization of the speed of convergence of   the trapezoidal rule, Numer. Math., 57(2), 123--138  (1990)
\bibitem{SCH72} G. Schmeisser, Optimale Quadraturformeln mit semidefiniten Kernen, Numer.  Math., 20  32--53 (1972/73)
\bibitem{SCH77} G. Schmeisser, A representation for the remainder of the Maclaurin quadrature  formula, Numer. Math., 27(3),  355--358 (1976/77)
\bibitem{SCH0} G. Schmeisser, Sampling, Gaussian quadrature, and Poisson summation formula, in ``Proceedings of the
1997 International Workshop on Sampling Theory and Applications'', Univ. Aveiro, Aveiro,  327--332 (1997)
\bibitem{SCH} G. Schmeisser, Quadrature over a semi-infinite interval and Mellin transform, in: Y.~Lyubarskii (ed.) ``Proceedings of the
1999 International Workshop on Sampling Theory and Applications'', 
(ISBN 82-7151-0991), Norwegian University of 
Science and Technology, Trondheim, 203--208 (1999)
\bibitem{S} F. Stenger,  Interpolation formulae based on the trapezoidal formula, J. Inst. Maths. Applics., 12, 103--114 (1973)
\bibitem{TIT} E.C. Titchmarsh, Introduction to the Theory of Fourier Integrals (2.~Ed.), Clarendon Press, Oxford (1948)
\end{thebibliography}


\end{document}